\documentclass{article}
\usepackage[numbers,sort&compress]{natbib}
\usepackage{amsmath}
\usepackage{amssymb}
\usepackage{multirow}
\usepackage{amsfonts}
\usepackage{mathrsfs}
\usepackage{booktabs}
\usepackage{authblk}
\usepackage{algorithmicx}
\usepackage{algorithm}
\usepackage{booktabs}
\usepackage{algpseudocode}
\usepackage{comment}
\usepackage{amsthm}
\usepackage{latexsym}
\usepackage{graphicx}
\usepackage{rotating}
\usepackage{subfig}
\usepackage{color}
\usepackage{varioref}
\usepackage{longtable}
\usepackage{epstopdf}
\allowdisplaybreaks[4]
\usepackage[colorlinks,linkcolor=blue]{hyperref}
\usepackage{cases}
\usepackage{subeqnarray}

\usepackage{geometry}
\begin{document} 

\title{A Novel Virus Diffusion Optimization (VDO) Algorithm for Global Optimization}
\author{Zhaoqi Sun\thanks{School of Mathematics and Computational Science, Xiangtan University, Xiangtan, 411105, China. Email: Angelia0606@163.com}, \hskip 0.2cm Qingsong Wang\thanks{Corresponding author. School of Mathematics and Computational Science, Xiangtan University, Xiangtan, 411105, China. Email: nothing2wang@hotmail.com}}
\date{\today}

\maketitle
	
\begin{abstract}
	Meta-heuristic algorithms are widely used to tackle complex optimization problems, including nonlinear, multimodal, and high-dimensional tasks. However, many existing methods suffer from premature convergence, limited exploration, and performance degradation in large-scale search spaces. To overcome these limitations, this paper introduces a novel Virus Diffusion Optimizer (VDO), inspired by the life-cycle and propagation dynamics of herpes-type viruses. VDO integrates four biologically motivated strategies, including viral tropism exploration, viral replication step regulation, virion diffusion propagation, and latency reactivation mechanism, to achieve a balanced trade-off between global exploration and local exploitation. Experiments on standard benchmark problems, including CEC 2017 and CEC 2022, demonstrate that VDO consistently surpasses state-of-the-art metaheuristics in terms of convergence speed, solution quality, and scalability. These results highlight the effectiveness of viral-inspired strategies in optimization and position VDO as a promising tool for addressing large-scale, complex problems in engineering and computational intelligence.To ensure reproducibility and foster further research, the source code of VDO is made publicly available.
\end{abstract}
	

\section{Introduction}
In recent years, optimization algorithms have been extensively applied across a wide range of fields, such as engineering design, industrial scheduling, healthcare, and data analytics. These methods have been proven effective in addressing complex, nonlinear, multimodal, and NP-hard optimization problems. For example, Absalom et al. \cite{EZUGWU2017189} proposed a hybrid method integrating simulated annealing with a symbiotic organisms search algorithm to solve the Traveling Salesman Problem (TSP). In another study, Pirozmand et al. \cite{10268392} decomposed the original Pathfinder Algorithm (PFA) into four sub-algorithms and developed the DFPA method for TSP through discretization and coupling techniques. Meanwhile, Fang Su et al. \cite{Su2025} introduced an Effective Dimension Extraction Mechanism (EDEM) that guides meta-heuristic algorithms using bidirectional mapping between high- and low-dimensional spaces. Experimental findings demonstrated that EDEM considerably enhanced the performance of various algorithms in tackling complex high-dimensional challenges and real-world applications.

Population-based stochastic search strategies, including Genetic Algorithms (GA) \cite{GA}, Particle Swarm Optimization (PSO) \cite{PSO}, and Artificial Bee Colony (ABC) \cite{ABC}, have been widely employed to approximate global optima without relying on gradient information. However, as established by the No-Free-Lunch theorem \cite{NFL}, no single algorithm is universally optimal for all types of problems. This underscores the importance of developing novel optimization strategies and hybrid approaches, particularly for challenging high-dimensional tasks. Metaheuristics are typically categorized according to their underlying sources of inspiration:
\begin{itemize}
	\item \textbf{Swarm Intelligence (SI)}—emulates collective intelligence observed in nature, such as bird flocking in Particle Swarm Optimization or bee foraging behavior in the Artificial Bee Colony algorithm.
	\item \textbf{Evolutionary Algorithms (EAs)}—simulate biological evolution through mechanisms including selection, crossover, and mutation, as seen in Genetic Algorithms and Differential Evolution.
	\item \textbf{Physics-Inspired Algorithms}—derive metaphors from physical laws or chemical processes, exemplified by Simulated Annealing, Gravitational Search Algorithm, and the Whale Optimization Algorithm.
	\item \textbf{Hybrid Methods}—combine elements from different metaheuristic paradigms to exploit their complementary strengths.
\end{itemize}

A fundamental challenge in designing metaheuristics lies in achieving an effective balance between exploration—broadly diversifying the search across the solution space—and exploitation—intensively refining promising regions—so as to avoid premature convergence or excessively slow convergence rates. As noted by Fazal et al. \cite{RN379}, the performance of a metaheuristic is largely governed by the interplay between intensification and diversification, making this equilibrium a pivotal aspect of algorithmic design.

Over the past decades, a variety of novel optimization algorithms have been proposed, many drawing inspiration from animal behaviors or natural phenomena \cite{Guo2023}. In parallel, significant enhancements have been made to existing metaheuristics to improve their performance and robustness. For instance, Houssein et al. \cite{HOUSSEIN2023119015} introduced a modified Sooty Tern Optimization Algorithm (mSTOA), which addressed the original method's limitations through an exploration–exploitation balancing strategy, adaptive parameter control, and a population reduction mechanism. In another study, Xiao et al. \cite{9953089} developed a Twin-Memory Bare-Bones Particle Swarm Optimization (TMBPSO) algorithm, incorporating a twin memory storage mechanism and a multiple memory retrieval strategy to empower the swarm with self-rectification capabilities, thereby alleviating premature convergence. Similarly, Li et al. \cite{GWO} proposed an Improved Gray Wolf Optimizer (IGWO), employing tent map-based initialization to enhance population diversity.

Despite these advances, high-dimensional optimization problems continue to present considerable difficulties. As the number of decision variables growing, the search space expands exponentially, making it difficult to maintain an appropriate balance between global exploration and local exploitation. In such contexts, many population-based algorithms tend to suffer from premature convergence or exhibit slow exploratory progress. Furthermore, complex engineering problems frequently incorporate constraints—such as nonlinear, multi-modal, or dynamic restrictions—that must be satisfied during the optimization of the objective function. The No-Free-Lunch theorem further reinforces the understanding that no single algorithm is capable of optimally solving every class of problem \cite{RN380,RN381}.

To overcome these limitations, this paper proposes the Virus Diffusion Optimizer (VDO), a novel metaheuristic algorithm inspired by the life-cycle and propagation behavior of herpes-type viruses. VDO integrates four biologically grounded strategies: viral tropism exploration, viral replication step regulation, virion diffusion propagation, and latency reactivation mechanism. These mechanisms work in concert to improve the trade-off between exploration and exploitation. In VDO, the population is categorized into different groups, each assigned a distinct update strategy, thereby enabling concurrent exploration across multiple directions in high-dimensional search spaces. Furthermore, the algorithm incorporates archive-driven reverse learning and adaptive step-size control to mitigate premature convergence and enhance solution quality.

Although traditional metaheuristics often perform well on low- to moderate-dimensional problems, their efficacy tends to decline significantly in high-dimensional settings due to the curse of dimensionality. This phenomenon undermines the discriminative capacity of distance measures and lowers the probability of sampling promising search directions. VDO alleviates these issues by sustaining population diversity and convergence efficiency even in problems involving hundreds to thousands of decision variables. As a result, it achieves robust performance in scenarios where conventional algorithms often fail.

The main contributions of this work are summarized as follows:

\begin{itemize}
	\item \textbf{A Novel Virus-Inspired Metaheuristic Framework.} We introduce a metaheuristic framework modeled after viral tropism, replication dynamics, and latency mechanisms to address complex optimization problems.
	
	\item \textbf{Detailed Algorithmic Mechanisms for Enhanced Search.} We design specific strategies—including adaptive step-size regulation, receptor-based particle selection, and archive-driven reverse learning—to improve search efficiency and convergence behavior.
	
	\item \textbf{Extensive Comparative Analysis and Performance Validation.} We conduct comprehensive comparisons with state-of-the-art metaheuristics, validating VDO's superior performance, particularly in high-dimensional optimization tasks.
\end{itemize}

The remainder of this paper is organized as follows. Section 2 describes the VDO algorithm and its biological foundations. Section 3 presents numerical experiments and comparative analyses. Finally, Section 4 concludes the paper and discusses potential directions for future research.

\section{Virus Diffusion Optimization}
High-dimensional and complex engineering problems typically involve vast search spaces, which demand substantial computational resources and pose challenges such as local optima and memory constraints. To address these difficulties, metaheuristic algorithms inspired by natural processes have been widely explored. Among them, the Virus Colony Search (VCS) algorithm simulates virus infection and diffusion strategies among host cells, providing an effective framework for exploring complex solution spaces \cite{LI201665}. 

Motivated by such biologically inspired strategies and drawing specifically on the lifecycle of Herpes Simplex Virus (HSV), we propose the Virus Diffusion Optimizer (VDO). VDO emulates key viral behaviors—including tropism, replication regulation, diffusion, and latency/reactivation—within a unified search framework, enabling adaptive and efficient navigation of high-dimensional optimization problems.

\subsection{Inspired by HSV}
Herpes simplex virus (HSV), a representative member of the \emph{Herpesviridae} family, offers a rich biological analogy for optimization. It is an enveloped double-stranded DNA virus whose glycoproteins recognize specific host-cell receptors, initiating infection with remarkable selectivity \cite{Kim2022}. After entry, HSV frequently establishes latency in sensory neurons and can reactivate under stress or immunosuppression, thereby sustaining long-term transmission \cite{Harrison2022,Fu2024}. 

The HSV lifecycle can be conceptually divided into four key stages: (1) receptor-mediated infection specificity, (2) regulated gene expression and replication cascades, (3) dual dissemination via budding and axonal transport, and (4) episodic latency followed by possible reactivation \cite{Feldman2002,Canova2024}. These strategies ensure both persistence and adaptability, making HSV a compelling metaphor for algorithmic design. In particular, they inspire mechanisms for targeted exploration, adaptive step-size control, hybrid local/global search, and diversity-preserving restarts in the proposed Virus Diffusion Optimizer (VDO). 

In VDO, each virion is modeled as a candidate solution in a high-dimensional search space, where fitness corresponds to replication efficiency as determined by the objective function. The search domain acts as the viral growth environment, and the algorithm's goal is to locate the optimal 'ecological niche' with the highest adaptability. By simulating HSV's infection, replication, propagation, and latency dynamics, VDO achieves an adaptive balance between exploration and exploitation, ensuring robust performance on large-scale, complex optimization tasks.

\subsection{Mathematical model and optimization algorithm}
In this section, inspired by the life cycle of herpes simplex virus (HSV), the Virus Diffusion Optimizer (VDO) is structured into four biologically motivated phases: viral tropism exploration, replication step regulation, virion diffusion propagation, and latency reactivation mechanism. Each phase models a distinct viral behavior, enabling the algorithm to adaptively navigate high-dimensional, multimodal landscapes.  

\subsubsection{Viral Tropism Exploration}
During HSV infection, the virions exhibit tropism, a selective affinity for specific host receptors such as nectin-1, HVEM, or 3-O-sulfated heparan sulfate, dictating which cells are preferentially infected. In VDO, this phase corresponds to the receptor-guided exploratory search immediately after population initialization. Candidate solutions selectively orient toward promising subregions of the search space through dimension-wise receptor selection controlled by \texttt{divide\_num}, while \texttt{tropism\_min} and \texttt{tropism\_max} modulate the probability of directional bias toward high-fitness areas.

The initial population is uniformly sampled within the search boundaries:
\begin{equation}
	\mathbf{x}_i^{(0)} = \mathbf{lb} + \text{rand}(1,d) \cdot (\mathbf{ub} - \mathbf{lb}), 
	\quad i = 1, \dots, N,
	\tag{1}
\end{equation}
where \(d\) is the dimension of the problem.

A subset of dimensions, representing receptor channels, is selected via the following.
\begin{equation}
	k = \Bigl\lfloor \frac{d}{\texttt{divide\_num}} \Bigr\rfloor, 
	\quad \mathcal{R} = \{r_1, \dots, r_k\},
	\tag{2}
\end{equation}
from a random permutation of all dimensions. Iterative filtering of candidate solutions is performed on the basis of the receptor dimension superiority.
\begin{equation}
	\mathcal{I}^{(j)} = \{\mathbf{x} \in \mathcal{I}^{(j-1)} : x[r_j] > x^*[r_j]\}, 
	\quad j = 1, \dots, k,
	\tag{3}
\end{equation}
with complementary selection applied if the surviving set shrinks excessively:
\begin{equation}
	\mathcal{I}^{(j)} = \mathcal{I}^{(j-1)} \setminus \mathcal{I}^{(j)},
	\quad \text{if } |\mathcal{I}^{(j)}| < 0.5 |\mathcal{I}^{(j-1)}|.
	\tag{4}
\end{equation}

Surviving particles (\(|\mathcal{I}^{(k)}| = M\)) are updated via receptor-guided steps:
\begin{equation}
	\mathbf{x}_i^{(t+1)} = \mathbf{x}_i^{(t)} + \eta \cdot S(\mathbf{x}_i^{(t)}, \mathbf{x}^*, \text{tropism}),
	\tag{5}
\end{equation}
where \(S(\cdot)\) implements receptor-based selective movement with intensity sampled from \([\texttt{tropism\_min}, \texttt{tropism\_max}]\).

This receptor-centric filtering and orientation mechanism
enables VDO to emulate HSV's biological tropism:
it selectively concentrates computational effort on subspaces with greater fitness potential while maintaining global diversity, 
thus laying the foundation for adaptive convergence in later phases.

\begin{figure}[H] 
	\centering 
	\includegraphics[width=0.5\textwidth]{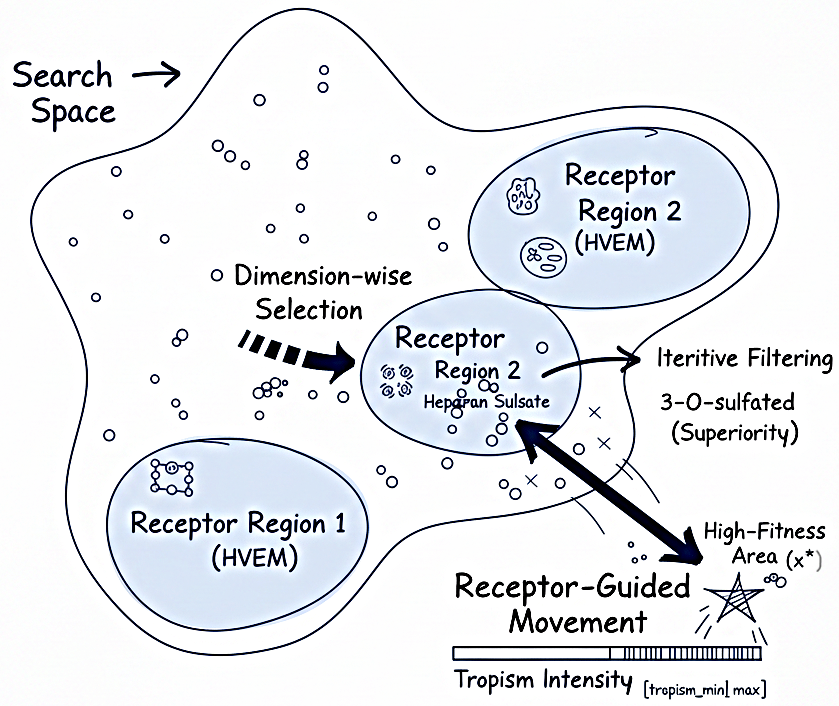} 
	\renewcommand{\figurename}{Figure}
	\caption{Viral Tropism Exploration.} 
	\label{Fig.Viral Tropism Exploration} 
\end{figure}

\subsubsection{Viral Replication Step Regulation}

When reactivated from latency, HSV enters a phase \emph{burst replication} inside the host cell, rapidly amplifying its genome and synthesizing structural proteins. 
This lytic process is characterized by an early stage of vigorous replication followed by gradual attenuation as host immune responses emerge. 
At the molecular level, the expression of the HSV gene follows a temporally ordered cascade: immediate–early (\(\alpha\)) genes activate early (\(\beta\)) genes that enable DNA replication and subsequently late (\(\gamma\)) genes direct the assembly of the capsid and envelope. 
This hierarchical timing of gene activation provides biological inspiration for a two-stage adaptive replication mechanism in the Virus Diffusion Optimizer (VDO).

In VDO, the \textbf{Burst Replication Phase} models the self-adaptive intensification of search agents around the best global solution, where the strength of the exploration decays over time in a biologically consistent manner. 
Let the current set of active agents after exploration of tropism be \(\mathcal{I} = \{\mathbf{x}_1, \dots, \mathbf{x}_M\}\), with global best \(\mathbf{x}^*\).
The mean environmental gradient that directs the collective movement toward the optimum is computed as:
\begin{equation}
	\Delta \mathbf{g} = \frac{1}{M} \sum_{\mathbf{x}_i \in \mathcal{I}} 
	\left( \mathbf{x}^* - \mathbf{x}_i \right).
	\tag{6}
\end{equation}
The overall replication strength is modulated by a time-decaying burst factor:
\begin{equation}
	\rho = w_0 * \beta_t,
	\tag{7}
\end{equation}
where \(w_0\) represents the initial replication strength.
Additionally, a secondary modulation term \(\beta_t\) refines this effect through nonlinear temporal scaling:
\begin{equation}
	\label{eq:burst_factor}
	\beta_t = \Biggl(1 - \frac{\mathrm{FEs}}{\mathrm{MaxFEs}}\Biggr)^{2\,\frac{\mathrm{FEs}}{\mathrm{MaxFEs}}},
	\tag{8}
\end{equation}
ensuring a smooth decay in replication intensity analogous to the declining viral load observed in HSV dynamics.

To emulate the two functional phases of viral gene expression, 
VDO employs both an \emph{exploratory “early” step} and an 
\emph{exploitative “late” step}.  
The exploratory update encourages population diversity:
\begin{equation}
	\mathbf{s}_{\mathrm{early}} 
	= \rho \, (\xi - 0.5) \, \beta_t \,
	\sin(2\pi \, \xi'),
	\tag{9}
\end{equation}
where \(\xi, \xi' \sim U(0,1)\) and \(\rho\) denotes the virion load parameter controlling replication amplitude. Conversely, the exploitative step focuses the population around promising regions:
\begin{equation}
	\mathbf{s}_{\mathrm{late}} 
	= 0.1 \, \rho \, (\xi - 0.5) \, \beta_t 
	\Biggl[ 1 + \frac{1}{2}\Bigl(1 + \tanh(R_0/\sqrt{1-R_0^2})\Bigr)\beta_t \Biggr],
	\tag{10}
\end{equation}
where \(R_0 = M/N\) denotes the effective reproduction ratio of active agents.
This adaptive modulation reflects the biological transition from early, high-rate viral genome replication to late, localized structural assembly. Each agent’s position is finally updated through a burst-oriented step that integrates envelope flipping and activation probability:
\begin{equation}
	\mathbf{x}_i^{(t+1)} 
	= \mathbf{x}_i^{(t)} 
	+ \rho \cdot 
	f\bigl(\Delta \mathbf{g}, \text{act\_prob}, \text{flip}\bigr),
	\tag{11}
\end{equation}
where \(f(\cdot)\) defines the stochastic step direction influenced by activation probability (\texttt{act\_prob}) and envelope flipping (\texttt{envelope\_flip}).
This formulation reproduces the collective and self-limiting replication dynamics of HSV: 
rapid expansion in early iterations followed by stabilized convergence as search pressure decays.

\begin{figure}[H] 
	\centering 
	\includegraphics[width=0.35\textwidth]{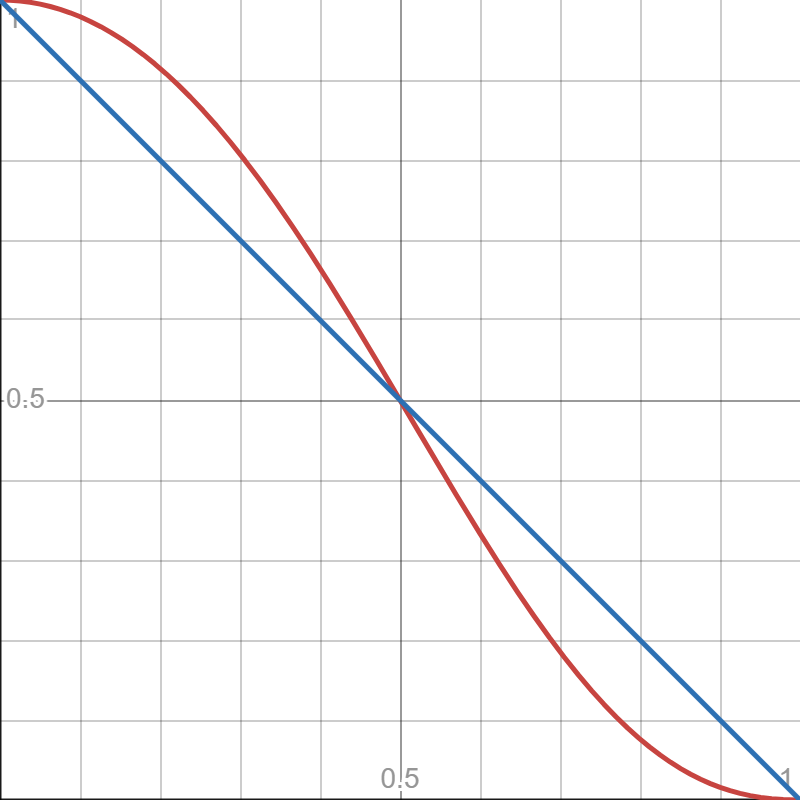} 
	\renewcommand{\figurename}{Figure}
	\caption{A comparison between the proposed adaptive method (red curve) and the conventional adaptive method (blue curve) for \(\beta_t\).} 
	\label{Fig.fun} 
\end{figure}

Equation~\eqref{eq:burst_factor} defines a time-dependent burst factor, \(\beta_t\), which adaptively modulates the step size of viral particles throughout the optimization process. At the early stage of optimization, when the number of function evaluations (\(\mathrm{FEs}\)) is small relative to the maximum allowed evaluations (\(\mathrm{MaxFEs}\)), \(\beta_t\) retains a relatively large value. This encourages broad exploration across the search space, thereby increasing the likelihood of discovering globally optimal regions. As the optimization progresses and \(\mathrm{FEs}\) approaches \(\mathrm{MaxFEs}\), \(\beta_t\) gradually decreases, reducing the step size and promoting fine-grained local search around high-quality solutions. Compared with traditional adaptive schemes that often rely on linear or monotonic decay, the nonlinear form of \(\beta_t\) in Equation~\eqref{eq:burst_factor} provides a more flexible and balanced control between exploration and exploitation, effectively enhancing convergence performance while mitigating the risk of premature stagnation. By alternating between \(\mathbf{s}_{\mathrm{early}}\) and \(\mathbf{s}_{\mathrm{late}}\), VDO emulates HSV's shift from rapid genome amplification to targeted structural assembly, adaptively balancing exploration and exploitation over the course of optimization.

\subsubsection{Virion Diffusion Propagation}

The \emph{egress phase} in HSV describes the release of mature virions from infected host cells through membrane budding or vesicle-mediated exocytosis, 
allowing them to infect neighboring cells and initiate secondary rounds of replication. 
This process involves two complementary dissemination pathways: 
(i) \emph{budding-mediated release}, which enables local spread, and 
(ii) \emph{cell-to-cell or axonal transmission}, which enables long-range propagation. 
Together, these pathways sustain infection diversity and spatial adaptability.

In the Virus Diffusion Optimizer (VDO), this biological phenomenon is abstracted as the \textbf{Viral Egress and Reinfection Phase}, 
corresponding to the algorithm’s global diffusion and diversity restoration mechanisms. 
After the replication burst, each particle \(\mathbf{x}_i \in \mathcal{I}\) undergoes one of two complementary update modes, reflecting HSV’s dual dissemination strategy.

\paragraph{(a) Budding-like update.}
With probability \(1 - p_{\mathrm{bud}}\), 
a particle performs a locally oriented propagation step:
\begin{equation}
	\mathbf{x}_i^{t+1} = 
	\mathbf{x}_i^{t} + \mathbf{s}_{\mathrm{step}} \odot \Delta \mathbf{g},
	\tag{12}
\end{equation}
where \(\mathbf{s}_{\mathrm{step}}\) is randomly selected from the exploratory (\(\mathbf{s}_{\mathrm{early}}\)) or exploitative (\(\mathbf{s}_{\mathrm{late}}\)) step definitions in the replication phase, 
\(\odot\) denotes element-wise multiplication, 
and \(\Delta \mathbf{g}\) is the global gradient computed in Eq.~(6).  
This operation simulates viral budding from the host membrane and local infection spread within proximal regions of the search space.

\paragraph{(b) Fusion or receptor-targeted jump.}
With complementary probability \(p_{\mathrm{bud}}\), 
a particle executes a receptor-based fusion or single-dimension jump:
\begin{equation}
	\text{either} \quad 
	\mathbf{x}_i[r] = x^*[r] + 0.1\,(\xi - 0.5)\,\beta_t\,\Delta g_r
	\quad \text{or} \quad 
	\mathbf{x}_i = (1 - \alpha)\,\mathbf{x}_i + \alpha\,\mathbf{x}^*,
	\tag{13}
\end{equation}
where \(r\) denotes a randomly selected receptor dimension, 
\(\xi \sim U(0,1)\), and \(\alpha \in \{0,1\}\) controls whether full or partial fusion occurs.  
This mechanism mirrors HSV’s cell-to-cell spread and receptor-driven targeting of new host environments.

\paragraph{(c) Multidirectional release and Lévy-based reinfection.}
VDO incorporates an \emph{envelope-flip} operator and \emph{Lévy-based long-distance jumps}, which aims to enhance global propagation.  
The envelope-flip operator (\texttt{envelope$\_$flip$\_$flags}) allows agents to reverse their movement direction, 
emulating multidirectional budding and enabling escape from local minima.  
Meanwhile, Lévy flights generate occasional long-range reinfection events:
\begin{equation}
	\mathbf{x}_i^{t+1} = 
	\mathbf{x}_i^{t} + L(d) \times (\mathbf{x}^* - \mathbf{x}_i^{t}),
	\tag{14}
\end{equation}
where \(L(d)\) denotes a Lévy-distributed random vector in \(d\)-dimensional space.
This stochastic reinfection mechanism introduces rare but significant global transitions, 
enhancing exploration and preventing premature convergence.

\paragraph{(d) Differential recombination and boundary reflection.}
In addition to diffusion, VDO integrates a Differential Evolution (DE)-based crossover, 
which functions analogously to viral genetic recombination during co-infection.  
This operator promotes structural diversity within the population.  
After each update, all particles are projected back into the feasible domain:
\begin{equation}
	\mathbf{x}_i^{t+1} =
	\min\!\bigl(\max(\mathbf{x}_i^{t+1}, \mathbf{lb}), \, \mathbf{ub}\bigr),
	\tag{15}
\end{equation}
ensuring biologically bounded yet dynamically adaptable exploration.

This phase thus completes the VDO infection cycle: 
it reintroduces diversity through budding, 
propagates high-fitness traits via fusion and recombination, 
and maintains long-range adaptability through Lévy reinfection. 
Collectively, these mechanisms reproduce the self-sustaining dynamics of HSV dissemination in multi-host environments, 
balancing local exploitation with global exploration.

\begin{figure}[H] 
	\centering 
	\includegraphics[width=0.7\textwidth]{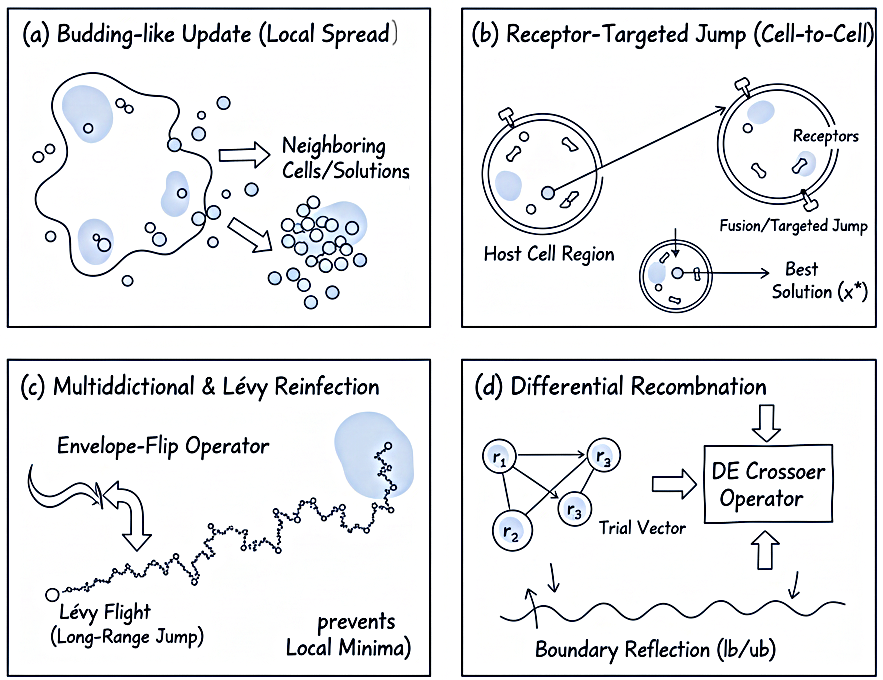} 
	\renewcommand{\figurename}{Figure}
	\caption{Virion Diffusion Propagation.} 
	\label{Fig.Virion Diffusion Propagation} 
\end{figure}

\subsubsection{Latency Reactivation Mechanism}

During the infection cycle of Herpes Simplex Virus (HSV), 
the virus can enter a \emph{latent state} within sensory neurons, 
where its genome persists as an episomal DNA in the host nucleus while remaining transcriptionally silent. 
This dormancy represents a delicate equilibrium between viral persistence and readiness for reactivation when environmental stimuli or immune suppression occur. 
The biological latency strategy ensures long-term survival, minimizes host immune detection, and preserves the viral population for future replication cycles.

In the Virus Diffusion Optimizer (VDO), this phenomenon is abstracted as the \textbf{Latency Maintenance and Reactivation Phase}, 
which governs algorithmic memory, stability, and controlled diversification. 
Here, a historical buffer $(rV, rV_{cos})$ stores several generations of solutions, forming a “latent reservoir’’ of promising search states. 
The parameter \texttt{latency\_depth} (denoted as $\tau_{\max}$) defines the retention span of this reservoir. 
Conceptually, this structure acts as a \emph{cryptobiotic archive} that preserves potentially valuable search patterns for later reactivation.

At each iteration, after fitness evaluation, the position and fitness of each particle $\mathbf{x}_i$ are recorded as:
\begin{equation}
	\mathcal{A}_i = 
	\{\mathbf{x}_i^{(t-\tau_{\max}+1)}, \dots, \mathbf{x}_i^{(t)}\}, 
	\quad
	\mathcal{F}_i = 
	\{f(\mathbf{x}_i^{(t-\tau_{\max}+1)}), \dots, f(\mathbf{x}_i^{(t)})\},
	\tag{16}
\end{equation}
where $\mathcal{A}_i$ and $\mathcal{F}_i$ respectively represent the archival positions and fitness trajectories of the $i$-th individual.  
This rolling buffer functions as a latent genomic memory, analogous to episomal persistence of HSV DNA.

When the latency counter $\tau$ exceeds its threshold $\tau_{\max}$, 
a reactivation event is triggered, reintroducing dormant solutions into the active population. 
The reactivation process proceeds as follows:
For each particle, the best historical position is identified:
\begin{equation}
	\mathbf{x}_i^{\mathrm{best}} = 
	\arg\min_{\mathbf{x}\in\mathcal{A}_i} f(\mathbf{x}),
	\tag{17}
\end{equation}
serving as a potential reactivation template analogous to the transcriptional reawakening of latent HSV genomes.
    
This latency–reactivation mechanism enables VDO to periodically restore lost diversity without external restart, 
preserving useful genetic material while preventing stagnation.  
By cycling between dormancy and renewal, 
the algorithm emulates HSV’s evolutionary strategy of alternating between latent persistence and lytic reactivation— 
achieving a dynamic balance between stability and adaptability in complex, multimodal search landscapes.

\begin{figure}[htbp] 
	\centering 
	\includegraphics[width=0.6\textwidth]{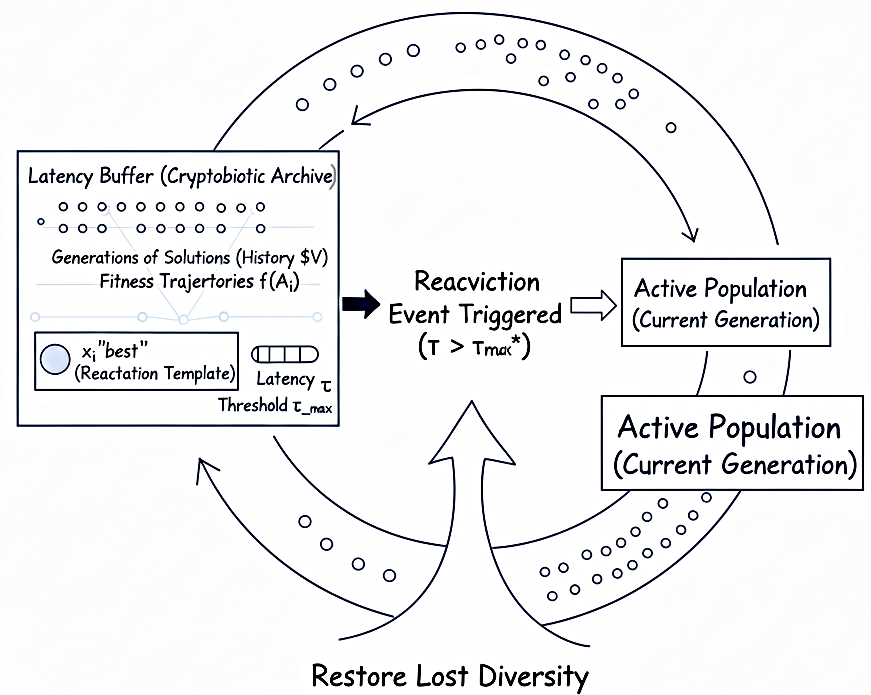} 
	\renewcommand{\figurename}{Figure}
	\caption{Latency Reactivation Mechanism.} 
	\label{Fig.Latency Reactivation Mechanism} 
\end{figure}

\subsection{Algorithm framework}
\begin{figure}[htbp] 
	\centering 
	\includegraphics[height=0.45\textwidth]{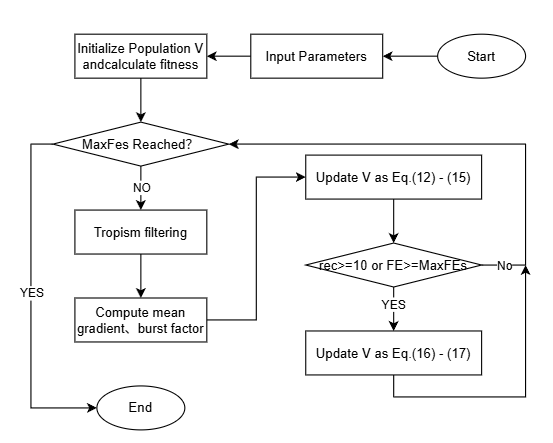} 
	\renewcommand{\figurename}{Figure}
	\caption{The overall flowchart of the Virus Diffusion Optimizer (VDO) algorithm.} 
	\label{Fig.process1} 
\end{figure}

Now we give the proposed algorithm as follows, see Figure \ref{Fig.process1} for the flowchart of the VDO algorithm and Algorithm \ref{alg:VDO_condensed} for details.  

\begin{algorithm}[H]
	\caption{Virus Diffusion Optimizer (VDO)}
	\label{alg:VDO_condensed}
	\begin{algorithmic}[1]
		\State \textbf{Inputs:}  $N$, $MaxFEs$, objective $f(\cdot)$.
		\State \textbf{Outputs:} best $\mathbf{x}^*$, $f^*$, Convergence\_curve
		\State Initialize $\{\mathbf{V}_i\}_{i=1}^N$ by Eq.(1); evaluate $f(V_i)$; set $\mathbf{x}^*, f^*, FEs\!:=\!N$, archive $rV, rV_{cost}$; $rec\!:=\!1$, $it\!:=\!1$
		
		\While{$FEs < MaxFEs$}
		\State \textbf{Tropism filtering:} obtain survivor index set $\mathcal{I}$ and $M=|\mathcal{I}|$  \Comment{Eq.(2)–(4)}
		\State Compute mean gradient $\Delta\mathbf{g}$ (Eq.(6)), time scale $\beta_t$ and $\text{burst\_factor}$ \Comment{Eq.(6)–(8)}
		
		\For{$i=1$ to $N$}
		\State $\mathbf{x}\leftarrow\mathbf{V}_i$
		\State \textbf{Burst update:} $\mathbf{x}\leftarrow \text{BurstUpdate}(\mathbf{x},\mathbf{x}^*,\Delta\mathbf{g},\text{burst\_factor},\beta_t)$ \Comment{Eq.(9)--(11)}
		\State \textbf{Diffusion:} optionally apply HSV’s dual dissemination strategy. \Comment{Eq.(12)--(15)}
		\State $f_{\text{new}}\leftarrow f(\mathbf{x})$; $FEs\!:=\!FEs+1$; store in archive $rV, rV_{cos}$ at slot $rec$
		\State If $f_{\text{new}}<f_i$ then $\mathbf{V}_i\leftarrow\mathbf{x}, f_i\leftarrow f_{\text{new}}$
		\State If $f_{\text{new}}<f^*$ then $\mathbf{x}^*\leftarrow\mathbf{x}, f^*\leftarrow f_{\text{new}}$
		\If{$FEs \ge MaxFEs$} 
		\State break 
		\EndIf
		\EndFor
		
		\State $rec\!:=\!rec+1$
		\If{$rec > \texttt{latency\_depth}$}
		
		\For{$i=1$ to $N$}
		\State $I^* \leftarrow \arg\min_j rV_{cos}(i,j)$ \Comment{Eq.(16)--(17)}
		\State $\mathbf{V}_i \leftarrow rV(i,:,I^*)$, $f_i \leftarrow rV_{cos}(i,I^*)$
		\EndFor
		
		\Return updated $\mathbf{V},\{f_i\}$
		\State $rec\!:=\!1$
		\EndIf
		
		\State $\text{Convergence\_curve}[it]\!:=\!f^*$; $it\!:=\!it+1$
		\EndWhile
		
		\Return $\mathbf{x}^*, f^*, \text{Convergence\_curve}$
	\end{algorithmic}
\end{algorithm}

\section{Experiments and Results}
In this section, we conduct numerical experiments using MATLAB on a PC equipped with an AMD Ryzen 9 7940H w (4.00 GHz), 16 GB of RAM, and the Windows 11 operating system. The experiments are designed to evaluate the performance of the Virus Diffusion Optimizer (VDO) on high-dimensional optimization problems using the CEC2017 and CEC2022 benchmark function sets. The source code, along with detailed documentation and examples to reproduce the experimental results in this paper, can be accessed at: \url{https://github.com/AngeliaSakura/VDO}. 

\subsection{Experimental method on CEC2017}
We choose CEC2017 \cite{CEC2017} function
test which contains 30 benchmark functions, including unimodal functions (F1–F2), simple multimodal functions (F3–F9), hybrid functions(F10–F19), and composition functions (F20–F30). Different varieties of functions can verify the comprehensive performance of the algorithm. Each function can be evaluated at dimensions of 10, 30, 50, and 100, with solving becoming increasingly difficult as the dimensionality increases. To validate the capability of VDO to address high-dimensional problems, we select the highest CEC2017 recommended dimension of 100 for experiments. In the experiment, a population size of 50 is utilized, with 30 independent run times conducted, each comprising 100,000 max number of function evolution.

\begin{table}[!ht]
	\centering
	\renewcommand{\tablename}{Table}
	\caption{Details of the CEC 2017.}
	\label{tab:cec2017}
	\begin{tabular}{llll}
		\toprule
		\textbf{ID} & \textbf{Function Equation} & \textbf{Class} & \textbf{Optimum} \\
		
		\midrule
		
		F1 & Shifted and Rotated Bent Cigar Function & Unimodal & 100 \\
		F2 & Shifted and Rotated Sum of Different Power Function & Unimodal & 200 \\
		F3 & Shifted and Rotated Zakharov Function & Unimodal & 300 \\

		\midrule
		
		F4 & Shifted and Rotated Rosenbrock’s Function & Multimodal & 400 \\
		F5 & Shifted and Rotated Rastrigin’s Function & Multimodal & 500 \\
		F6 & Shifted and Rotated Expanded Scaffer’s F6 Function & Multimodal & 600 \\
		F7 & Shifted and Rotated Lunacek Bi-Rastrigin Function & Multimodal & 700 \\
		F8 & Shifted and Rotated Non-Continuous Rastrigin’s Function & Multimodal & 800 \\
		F9 & Shifted and Rotated Lévy Function & Multimodal & 900 \\
		F10 & Shifted and Rotated Schwefel’s Function & Multimodal & 1000 \\

		\midrule
		
		F11 & Hybrid Function 1 (N = 3) & Hybrid & 1100 \\
		F12 & Hybrid Function 2 (N = 3) & Hybrid & 1200 \\
		F13 & Hybrid Function 3 (N = 3) & Hybrid & 1300 \\
		F14 & Hybrid Function 4 (N = 4) & Hybrid & 1400 \\
		F15 & Hybrid Function 5 (N = 4) & Hybrid & 1500 \\
		F16 & Hybrid Function 6 (N = 4) & Hybrid & 1600 \\
		F17 & Hybrid Function 6 (N = 5) & Hybrid & 1700 \\
		F18 & Hybrid Function 6 (N = 5) & Hybrid & 1800 \\
		F19 & Hybrid Function 6 (N = 5) & Hybrid & 1900 \\
		F20 & Hybrid Function 6 (N = 6) & Hybrid & 2000 \\

		\midrule
		
		F21 & Composition Function 1 (N = 3) & Composition & 2100 \\
		F22 & Composition Function 2 (N = 3) & Composition & 2200 \\
		F23 & Composition Function 3 (N = 4) & Composition & 2300 \\
		F24 & Composition Function 4 (N = 4) & Composition & 2400 \\
		F25 & Composition Function 5 (N = 5) & Composition & 2500 \\
		F26 & Composition Function 6 (N = 5) & Composition & 2600 \\
		F27 & Composition Function 7 (N = 6) & Composition & 2700 \\
		F28 & Composition Function 8 (N = 6) & Composition & 2800 \\
		F29 & Composition Function 9 (N = 3) & Composition & 2900 \\
		F30 & Composition Function 10 (N = 3) & Composition & 3000 \\
		\bottomrule
	\end{tabular}
\end{table}

Eleven optimization algorithms are selected as control groups.
To ensure fairness in the comparison, the parameters for all selected algorithms are set according to Table \ref{tab:control_group_params_refined}.

\begin{table}[!hb]
	\centering
	\renewcommand{\tablename}{Table}
	\caption{Key parameters and original year of control-group algorithms.}
	\label{tab:control_group_params_refined}
	\begin{tabular}{llll}
		\toprule
		\textbf{Algorithm} 
		& \textbf{Key parameters / settings} 
		& \textbf{Original year} & \textbf{Reference} \\
		\midrule
		GA        & \( p_c = 0.8,\; p_m = 0.05 \) 
		& 1976 & \cite{GA}\\
		PSO       & \( c_1 = 2,\; c_2 = 2 \) 
		& 1995 & \cite{PSO} \\
		GWO       & \( a \in [2,\,0] \) 
		& 2014 & \cite{GWO} \\
		WOA       & \( a_1 \in [2,\,0],\; a_2 \in [-2, -1],\; b = 1 \) 
		& 2016 & \cite{WOA} \\
		HHO       & \( E_0 \in [-1,1] \)  
		& 2019 & \cite{HHO} \\
		HO        &  -
		& 2024 & \cite{HO} \\
		NRBO      & \( DF = 0.6 \) 
		& 2024 & \cite{NRBO} \\
		FATA      & \( arf = 0.2 \) 
		& 2024 & \cite{FATA} \\
		MGO       & \( w = 2,\; rec\_num = 10 \) 
		& 2024 & \cite{MGO} \\
		LEA       & \( h_{\max} = 0.7,\; h_{\min} = 0,\; \lambda_c = 0.5,\; \lambda_p = 0.5 \) 
		& 2024 & \cite{LEA} \\
		ALA       & \( \mathrm{vec\_flag} = [1, -1] \) 
		& 2025 & \cite{ALA} \\
		\bottomrule
	\end{tabular}
\end{table}

Based on the experimental results, we evaluate the algorithms' performance using the mean and variance of their optimization outcomes. The mean value reflects the average solution quality achieved by each algorithm, while the variance indicates the stability of its performance across multiple runs. These metrics provide comprehensive insights into both the accuracy and reliability of the optimization methods under comparison.

\begin{sidewaystable}[htbp]
	\centering
	\renewcommand{\tablename}{Table}
	\caption{Ablation experiment results of VDO on F1–F10 in the CEC2017.}
	\label{tab:cec2017_F1_F10}
	{\scriptsize
		\begin{tabular}{llllllllllllll}
			\toprule
			\textbf{F} & \textbf{Data type} & \textbf{VDO} & \textbf{GA} & \textbf{PSO} & \textbf{GWO} & \textbf{WOA} & \textbf{HHO} & \textbf{HO} & \textbf{NRBO} & \textbf{FATA} & \textbf{MGO} & \textbf{LEA} & \textbf{ALA} \\
			\midrule
			
			\( F_1 \) & Mean & 3.82E+06 & 1.30E+09 & 3.97E+10 & 4.34E+10 & 1.07E+10 & 6.61E+08 & 1.07E+10 & 1.49E+11 & 3.69E+10 & 2.11E+08 & 1.31E+11 & 1.84E+07 \\
			& Variance & 3.03E+12 & 3.93E+16 & 1.74E+20 & 8.33E+19 & 6.15E+18 & 6.33E+15 & 7.01E+18 & 1.61E+20 & 5.47E+19 & 1.93E+16 & 3.60E+20 & 2.59E+14 \\
			& Rank\_M & 1 & 5 & 9 & 10 & 7 & 4 & 6 & 12 & 8 & 3 & 11 & 2 \\
			& Rank\_V & 1 & 5 & 11 & 9 & 6 & 3 & 7 & 10 & 8 & 4 & 12 & 2 \\
			
			\( F_2 \) & Mean & 4.76E+74 & 1.31E+120 & 1.24E+130 & 1.05E+157 & 2.17E+113 & 3.54E+138 & 6.28E+144 & 3.52E+139 & 1.70E+132 & 6.40E+165 & 7.79E+98 & 1.92E+99 \\
			& Variance & 1.53E+150 & 1.21E+241 & 1.08E+261 & Inf & 1.53E+227 & 7.29E+277 & 2.17E+290 & 5.08E+279 & 1.97E+265 & Inf & 4.25E+198 & 1.51E+199 \\
			& Rank\_M & 1 & 5 & 6 & 11 & 4 & 8 & 10 & 9 & 7 & 12 & 2 & 3 \\
			& Rank\_V & 1 & 5 & 6 & 11 & 4 & 8 & 10 & 9 & 7 & 11 & 2 & 3 \\
			
			\( F_3 \) & Mean & 4.60E+05 & 1.51E+05 & 4.62E+05 & 2.62E+05 & 9.42E+05 & 2.27E+05 & 2.76E+05 & 2.92E+05 & 2.70E+05 & 7.00E+05 & 1.01E+06 & 1.52E+05 \\
			& Variance & 2.82E+09 & 1.05E+09 & 1.34E+10 & 7.34E+08 & 2.02E+10 & 4.41E+08 & 1.90E+08 & 7.52E+08 & 5.54E+08 & 4.29E+09 & 4.06E+10 & 4.23E+08 \\
			& Rank\_M & 8 & 1 & 9 & 4 & 11 & 3 & 6 & 7 & 5 & 10 & 12 & 2 \\
			& Rank\_V & 8 & 7 & 10 & 5 & 11 & 3 & 1 & 6 & 4 & 9 & 12 & 2 \\
			
			\( F_4 \) & Mean & 8.41E+02 & 1.11E+03 & 3.94E+03 & 3.56E+03 & 3.22E+03 & 1.08E+03 & 3.47E+03 & 2.19E+04 & 3.91E+03 & 8.09E+02 & 1.87E+04 & 8.56E+02 \\
			& Variance & 2.14E+03 & 9.07E+03 & 2.46E+06 & 8.11E+05 & 3.77E+05 & 1.02E+04 & 6.03E+05 & 8.29E+06 & 7.49E+05 & 4.29E+03 & 1.46E+07 & 1.03E+04 \\
			& Rank\_M & 2 & 5 & 10 & 8 & 6 & 4 & 7 & 12 & 9 & 1 & 11 & 3 \\
			& Rank\_V & 1 & 3 & 10 & 9 & 6 & 4 & 7 & 11 & 8 & 2 & 12 & 5 \\
			
			\( F_5 \) & Mean & 9.41E+02 & 1.71E+03 & 1.03E+03 & 1.11E+03 & 1.64E+03 & 1.55E+03 & 1.35E+03 & 1.96E+03 & 1.63E+03 & 1.32E+03 & 2.37E+03 & 1.02E+03 \\
			& Variance & 3.87E+03 & 1.39E+04 & 6.16E+03 & 6.23E+03 & 1.02E+04 & 5.55E+03 & 1.89E+03 & 8.00E+03 & 1.67E+03 & 2.18E+03 & 1.59E+04 & 6.75E+03 \\
			& Rank\_M & 1 & 10 & 3 & 4 & 9 & 7 & 6 & 11 & 8 & 5 & 12 & 2 \\
			& Rank\_V & 4 & 11 & 6 & 7 & 10 & 5 & 2 & 9 & 1 & 3 & 12 & 8 \\
			
			\( F_6 \) & Mean & 6.09E+02 & 6.84E+02 & 6.33E+02 & 6.36E+02 & 6.98E+02 & 6.83E+02 & 6.68E+02 & 7.02E+02 & 6.80E+02 & 6.09E+02 & 7.19E+02 & 6.31E+02 \\
			& Variance & 1.08E+00 & 6.67E+01 & 1.13E+02 & 8.20E+00 & 7.69E+01 & 2.28E+01 & 2.14E+01 & 2.85E+01 & 1.22E+01 & 1.18E+00 & 4.12E+01 & 3.82E+01 \\
			& Rank\_M & 1 & 9 & 4 & 5 & 10 & 8 & 6 & 11 & 7 & 2 & 12 & 3 \\
			& Rank\_V & 1 & 10 & 12 & 3 & 11 & 6 & 5 & 7 & 4 & 2 & 9 & 8 \\
			
			\( F_7 \) & Mean & 1.38E+03 & 2.22E+03 & 1.61E+03 & 1.97E+03 & 3.55E+03 & 3.60E+03 & 3.04E+03 & 3.61E+03 & 3.12E+03 & 1.78E+03 & 5.79E+03 & 1.70E+03 \\
			& Variance & 9.81E+03 & 4.38E+04 & 4.09E+04 & 1.99E+04 & 1.81E+04 & 2.46E+04 & 1.75E+04 & 3.87E+04 & 3.32E+04 & 3.58E+03 & 1.49E+05 & 4.05E+04 \\
			& Rank\_M & 1 & 6 & 2 & 5 & 9 & 10 & 7 & 11 & 8 & 4 & 12 & 3 \\
			& Rank\_V & 2 & 11 & 10 & 5 & 4 & 6 & 3 & 8 & 7 & 1 & 12 & 9 \\
			
			\( F_8 \) & Mean & 1.19E+03 & 2.07E+03 & 1.35E+03 & 1.41E+03 & 2.08E+03 & 1.97E+03 & 1.78E+03 & 2.39E+03 & 2.08E+03 & 1.60E+03 & 2.68E+03 & 1.26E+03 \\
			& Variance & 3.87E+03 & 2.37E+04 & 1.41E+04 & 4.31E+03 & 7.19E+03 & 3.67E+03 & 6.48E+03 & 1.07E+04 & 2.59E+03 & 4.18E+03 & 5.88E+04 & 1.08E+04 \\
			& Rank\_M & 1 & 8 & 3 & 4 & 10 & 7 & 6 & 11 & 9 & 5 & 12 & 2 \\
			& Rank\_V & 3 & 11 & 10 & 5 & 7 & 2 & 6 & 8 & 1 & 4 & 12 & 9 \\
			
			\( F_9 \) & Mean & 5.99E+03 & 6.02E+04 & 4.57E+04 & 3.36E+04 & 5.49E+04 & 5.16E+04 & 2.83E+04 & 6.42E+04 & 5.54E+04 & 2.60E+04 & 1.65E+05 & 1.76E+04 \\
			& Variance & 1.40E+06 & 3.41E+08 & 6.67E+08 & 1.27E+08 & 2.10E+08 & 2.65E+07 & 6.56E+06 & 1.42E+07 & 2.97E+07 & 2.63E+07 & 1.35E+08 & 1.41E+07 \\
			& Rank\_M & 1 & 10 & 6 & 5 & 8 & 7 & 4 & 11 & 9 & 3 & 12 & 2 \\
			& Rank\_V & 1 & 11 & 12 & 8 & 10 & 6 & 2 & 4 & 7 & 5 & 9 & 3 \\
			
			\( F_{10} \) & Mean & 1.60E+04 & 2.25E+04 & 1.57E+04 & 1.62E+04 & 2.41E+04 & 2.14E+04 & 1.79E+04 & 2.98E+04 & 3.13E+04 & 2.46E+04 & 3.01E+04 & 1.84E+04 \\
			& Variance & 1.27E+06 & 4.98E+06 & 2.12E+06 & 4.67E+06 & 5.47E+06 & 2.60E+06 & 4.38E+06 & 1.64E+06 & 5.58E+05 & 1.27E+06 & 2.80E+06 & 5.90E+06 \\
			& Rank\_M & 2 & 7 & 1 & 3 & 8 & 6 & 4 & 10 & 12 & 9 & 11 & 5 \\
			& Rank\_V & 2 & 10 & 5 & 9 & 11 & 6 & 8 & 4 & 1 & 3 & 7 & 12 \\
			
			\bottomrule
		\end{tabular}
	}
\end{sidewaystable}

\begin{figure}[H]
\setlength\tabcolsep{0.3pt}
\centering
\begin{tabular}{c}
\includegraphics[width=0.98\textwidth]{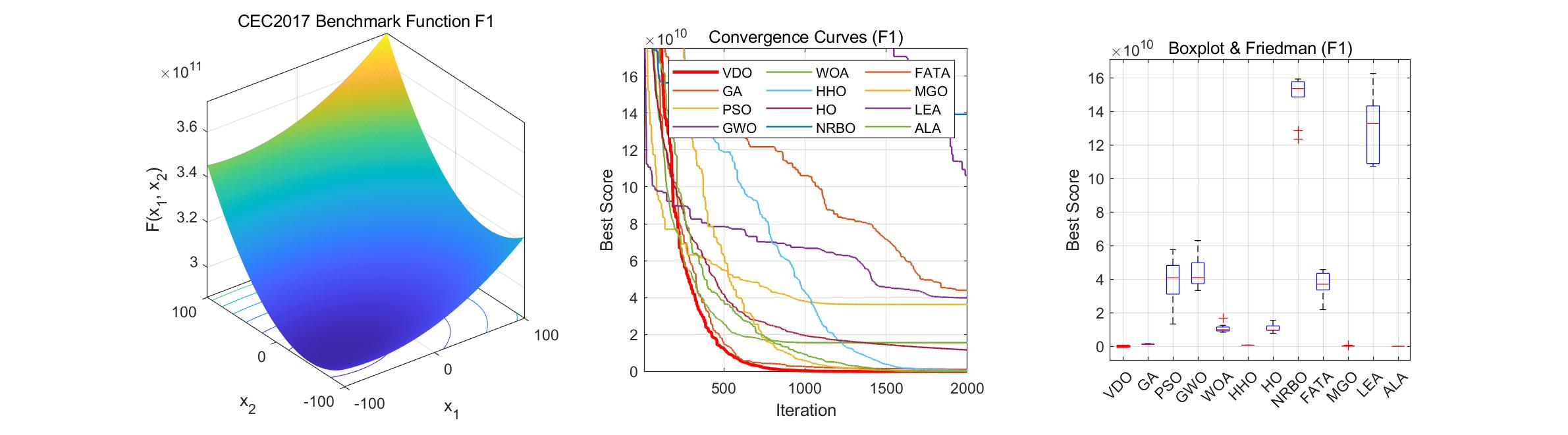}\\
\includegraphics[width=0.98\linewidth]{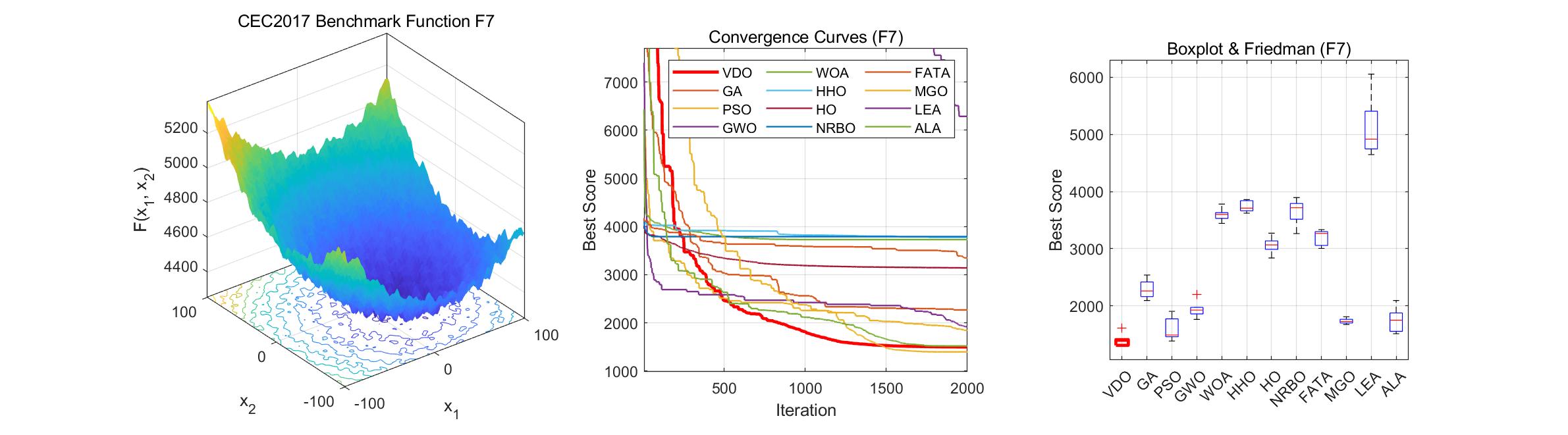}\\
\includegraphics[width=0.98\linewidth]{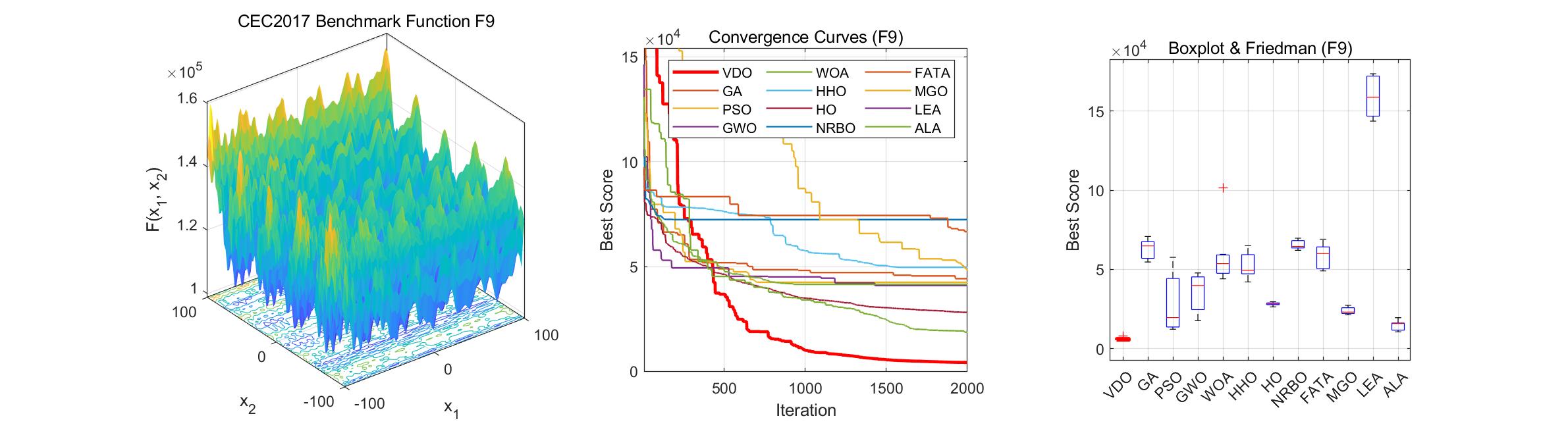}
\end{tabular}
\caption{Convergence curves of VDO and other algorithms on CEC2017 benchmark functions F1, F7, and F9, respectively.}
\label{fig:F1F7F9}
\end{figure}

	
	
	
	
	
	

\begin{figure}[!ht]
\setlength\tabcolsep{0.3pt}
\centering
\begin{tabular}{c}
\includegraphics[width=0.98\textwidth]{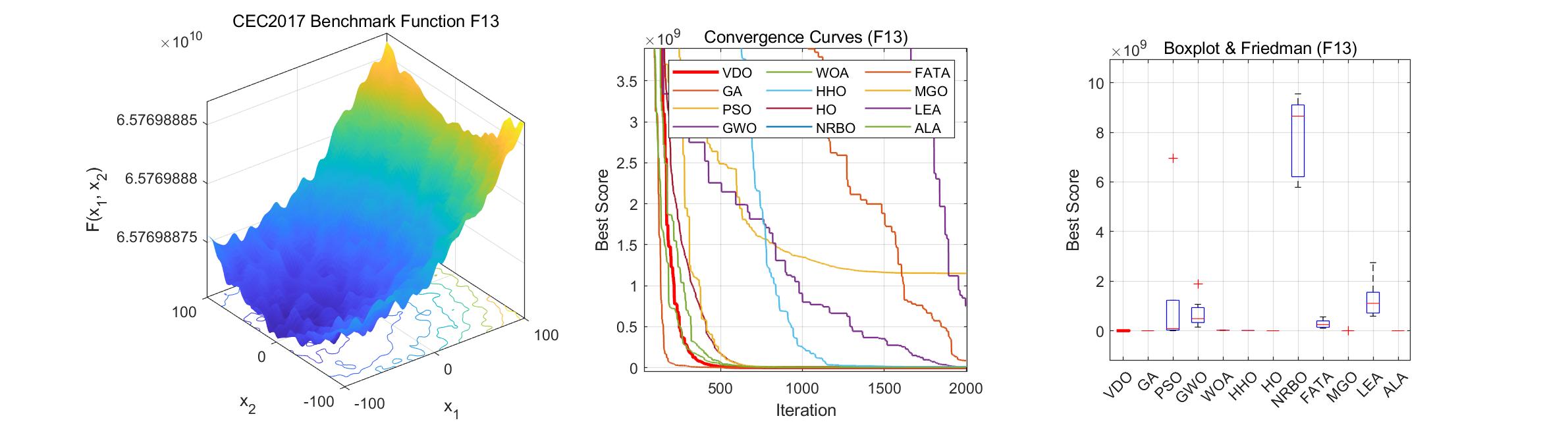}\\
\includegraphics[width=0.98\linewidth]{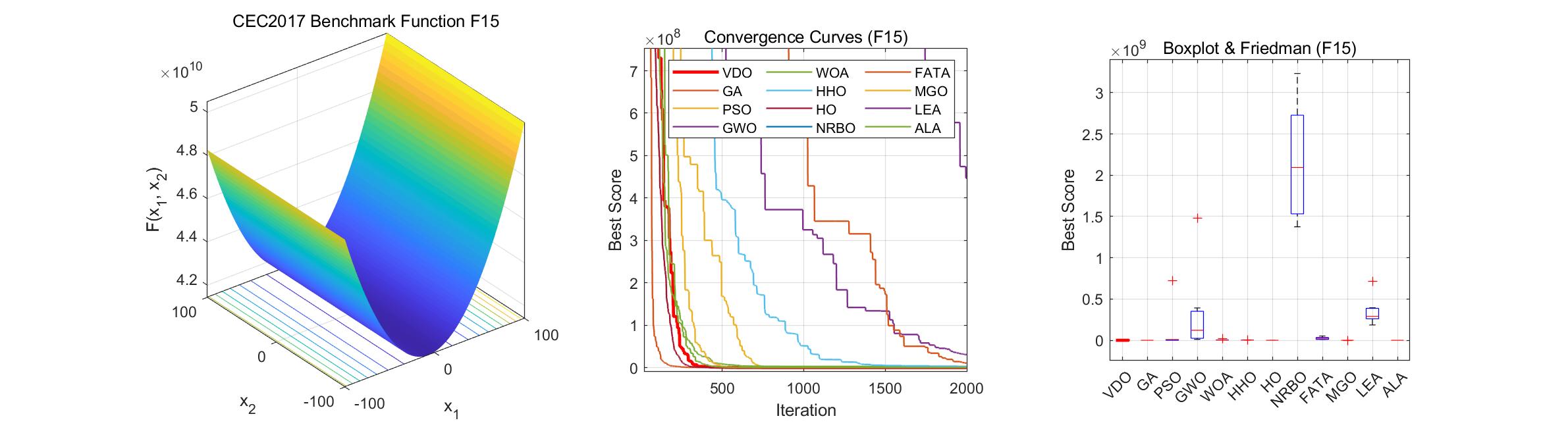}\\
\includegraphics[width=0.98\linewidth]{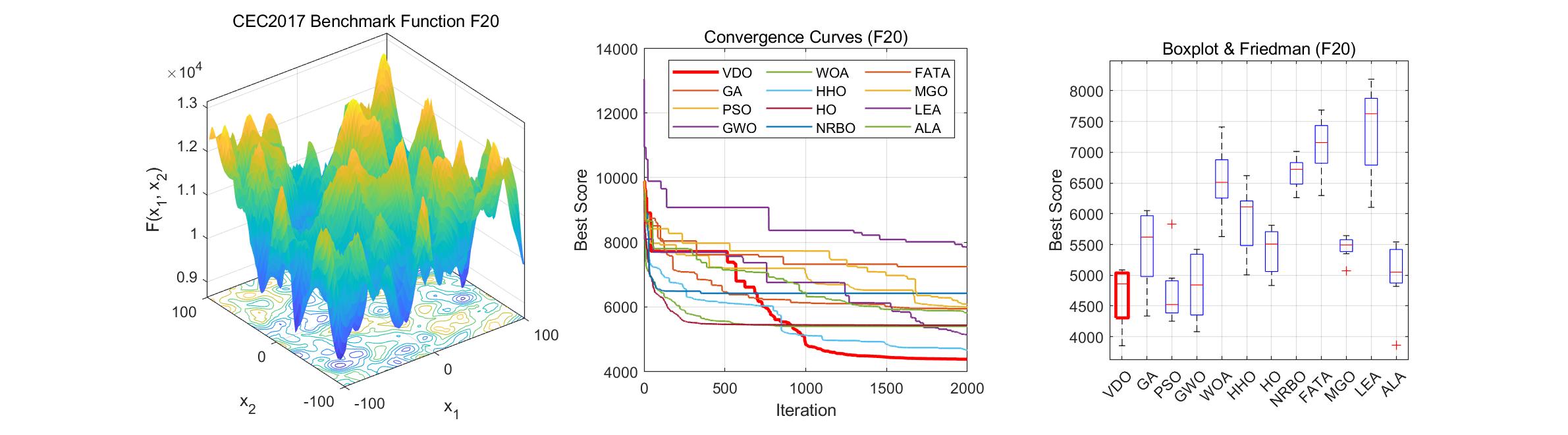}
\end{tabular}
\caption{Convergence curves of VDO and other algorithms on CEC2017 benchmark functions F13, F15, and F20.}
\label{fig:F13F15F20}
\end{figure}

\begin{figure}[!ht]
\setlength\tabcolsep{0.3pt}
\centering
\begin{tabular}{c}
\includegraphics[width=0.98\textwidth]{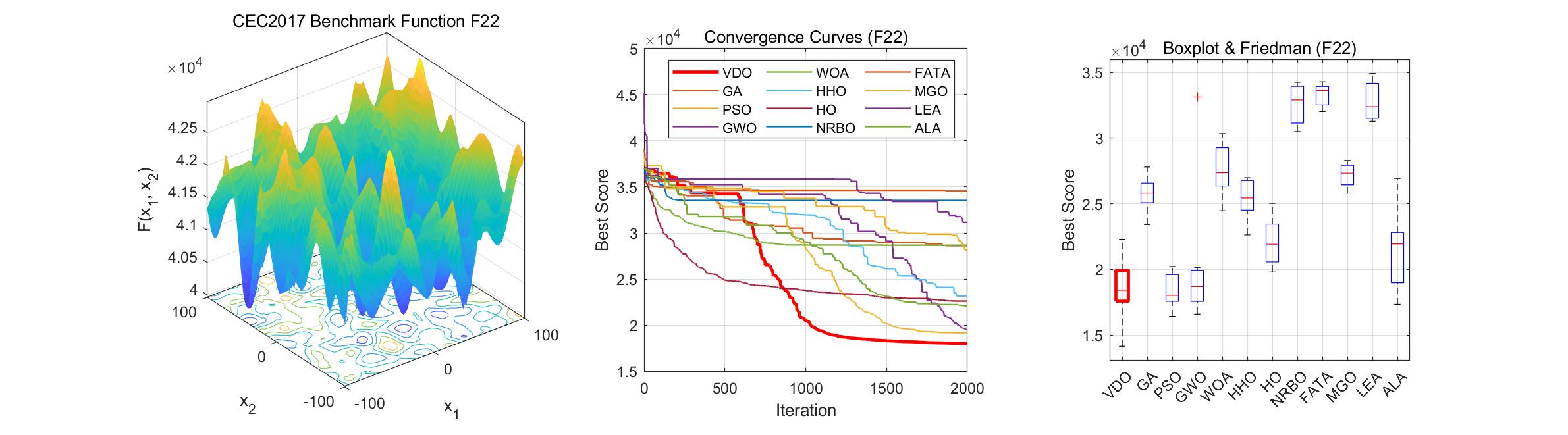}\\
\includegraphics[width=0.98\linewidth]{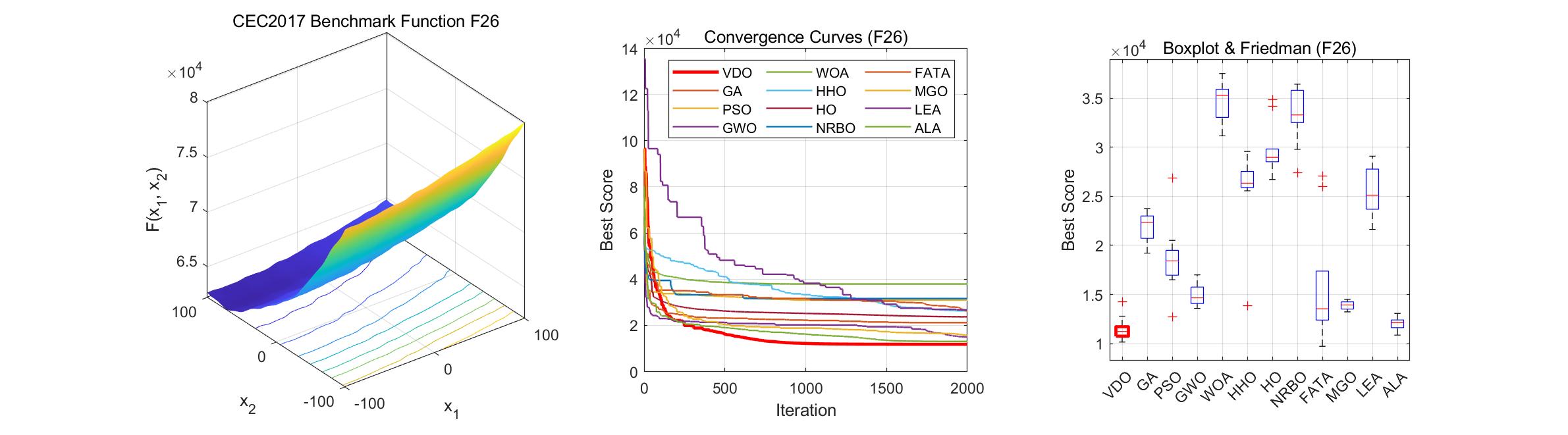}\\
\includegraphics[width=0.98\linewidth]{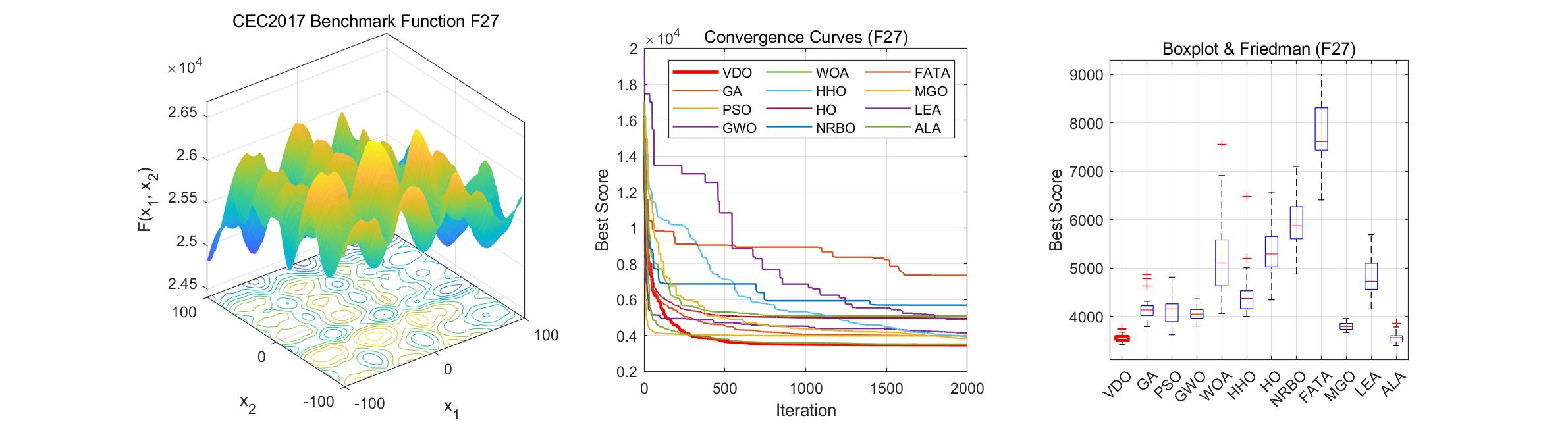}
\end{tabular}
\caption{Convergence curves of VDO and other algorithms on CEC2017 benchmark functions F22, F26, and F27.}
\label{fig:F7-F9}
\end{figure}

\begin{sidewaystable}[htbp]
	\centering
	\renewcommand{\tablename}{Table}
	\caption{Ablation experiment results of VDO on F11–F20 in the CEC2017.}
	\label{tab:cec2017_F11_F20}
	{\scriptsize
		\begin{tabular}{llllllllllllll}
			\toprule
			\textbf{F} & \textbf{Data type} & \textbf{VDO} & \textbf{GA} & \textbf{PSO} & \textbf{GWO} & \textbf{WOA} & \textbf{HHO} & \textbf{HO} & \textbf{NRBO} & \textbf{FATA} & \textbf{MGO} & \textbf{LEA} & \textbf{ALA} \\
			\midrule
			
			\( F_{11} \) & Mean & 2.74E+04 & 7.76E+03 & 9.99E+03 & 6.38E+04 & 1.28E+05 & 1.12E+04 & 5.41E+04 & 8.32E+04 & 4.98E+04 & 1.07E+05 & 3.79E+05 & 4.91E+03 \\
			& Variance & 2.92E+07 & 1.16E+07 & 1.22E+08 & 5.34E+07 & 3.75E+09 & 4.04E+07 & 5.86E+07 & 3.01E+07 & 4.93E+08 & 3.93E+08 & 3.54E+09 & 8.05E+05 \\
			& Rank\_M & 5 & 2 & 3 & 8 & 11 & 4 & 7 & 9 & 6 & 10 & 12 & 1 \\
			& Rank\_V & 3 & 2 & 8 & 6 & 12 & 5 & 7 & 4 & 10 & 9 & 11 & 1 \\
			
			\( F_{12} \) & Mean & 5.51E+07 & 4.24E+08 & 1.05E+10 & 9.15E+09 & 2.87E+09 & 7.15E+08 & 1.96E+09 & 5.70E+10 & 7.36E+09 & 3.36E+08 & 1.98E+10 & 3.24E+07 \\
			& Variance & 3.63E+14 & 1.28E+16 & 4.73E+19 & 2.69E+19 & 8.64E+17 & 2.33E+16 & 6.15E+17 & 1.21E+20 & 1.11E+19 & 6.83E+15 & 2.32E+19 & 1.82E+14 \\
			& Rank\_M & 2 & 4 & 10 & 9 & 7 & 5 & 6 & 12 & 8 & 3 & 11 & 1 \\
			& Rank\_V & 2 & 4 & 11 & 10 & 7 & 5 & 6 & 12 & 8 & 3 & 9 & 1 \\
			
			\( F_{13} \) & Mean & 9.08E+03 & 1.76E+04 & 1.37E+09 & 6.91E+08 & 1.87E+07 & 9.17E+06 & 7.10E+04 & 7.90E+09 & 2.77E+08 & 2.31E+04 & 1.26E+09 & 1.39E+04 \\
			& Variance & 3.88E+07 & 9.55E+07 & 6.36E+18 & 3.60E+17 & 3.98E+13 & 3.92E+12 & 8.11E+08 & 2.52E+18 & 2.81E+16 & 1.10E+09 & 5.60E+17 & 6.54E+07 \\
			& Rank\_M & 1 & 3 & 11 & 9 & 7 & 6 & 5 & 12 & 8 & 4 & 10 & 2 \\
			& Rank\_V & 1 & 3 & 12 & 9 & 7 & 6 & 4 & 11 & 8 & 5 & 10 & 2 \\
			
			\( F_{14} \) & Mean & 9.26E+05 & 2.53E+06 & 2.39E+06 & 6.11E+06 & 6.95E+06 & 2.70E+06 & 3.08E+06 & 1.39E+07 & 4.67E+06 & 6.46E+06 & 3.03E+07 & 8.01E+03 \\
			& Variance & 2.97E+11 & 1.15E+12 & 2.99E+12 & 1.08E+13 & 8.38E+12 & 6.68E+11 & 1.69E+12 & 2.74E+13 & 3.02E+12 & 7.03E+12 & 2.37E+14 & 3.95E+07 \\
			& Rank\_M & 2 & 4 & 3 & 8 & 10 & 5 & 6 & 11 & 7 & 9 & 12 & 1 \\
			& Rank\_V & 2 & 4 & 6 & 10 & 9 & 3 & 5 & 11 & 7 & 8 & 12 & 1 \\
			
			\( F_{15} \) & Mean & 4.51E+03 & 9.39E+03 & 5.59E+08 & 1.36E+08 & 4.56E+06 & 2.21E+06 & 7.42E+04 & 3.19E+09 & 1.53E+07 & 2.95E+04 & 3.26E+08 & 9.30E+03 \\
			& Variance & 8.44E+06 & 3.96E+07 & 5.08E+17 & 1.04E+16 & 3.01E+13 & 2.47E+11 & 8.00E+08 & 1.40E+18 & 6.36E+13 & 4.06E+09 & 2.61E+16 & 8.76E+07 \\
			& Rank\_M & 1 & 3 & 11 & 9 & 7 & 6 & 5 & 12 & 8 & 4 & 10 & 2 \\
			& Rank\_V & 1 & 2 & 11 & 9 & 7 & 6 & 4 & 12 & 8 & 5 & 10 & 3 \\
			
			\( F_{16} \) & Mean & 5.54E+03 & 6.43E+03 & 5.95E+03 & 5.93E+03 & 1.35E+04 & 7.71E+03 & 9.82E+03 & 1.37E+04 & 7.45E+03 & 7.31E+03 & 1.23E+04 & 5.91E+03 \\
			& Variance & 2.46E+05 & 2.12E+06 & 5.87E+05 & 4.65E+05 & 2.35E+06 & 1.65E+06 & 1.48E+06 & 5.00E+06 & 4.38E+05 & 6.25E+05 & 1.59E+06 & 4.24E+05 \\
			& Rank\_M & 1 & 5 & 4 & 3 & 11 & 8 & 9 & 12 & 7 & 6 & 10 & 2 \\
			& Rank\_V & 1 & 10 & 5 & 4 & 11 & 9 & 7 & 12 & 3 & 6 & 8 & 2 \\
			
			\( F_{17} \) & Mean & 4.69E+03 & 5.82E+03 & 6.19E+03 & 4.91E+03 & 8.54E+03 & 6.23E+03 & 7.00E+03 & 1.79E+04 & 6.96E+03 & 5.34E+03 & 9.06E+03 & 4.98E+03 \\
			& Variance & 3.11E+05 & 3.33E+05 & 2.41E+06 & 5.64E+05 & 6.59E+05 & 8.99E+05 & 2.24E+05 & 1.92E+08 & 2.06E+06 & 2.09E+05 & 4.18E+05 & 3.56E+05 \\
			& Rank\_M & 1 & 5 & 6 & 2 & 10 & 7 & 9 & 12 & 8 & 4 & 11 & 3 \\
			& Rank\_V & 3 & 4 & 11 & 7 & 8 & 9 & 2 & 12 & 10 & 1 & 6 & 5 \\
			
			\( F_{18} \) & Mean & 1.80E+06 & 5.21E+06 & 3.67E+06 & 6.02E+06 & 6.25E+06 & 3.85E+06 & 2.59E+06 & 1.99E+07 & 6.54E+06 & 9.72E+06 & 5.18E+07 & 1.12E+05 \\
			& Variance & 1.42E+12 & 5.95E+12 & 4.35E+12 & 1.09E+13 & 1.19E+13 & 5.18E+12 & 1.59E+12 & 1.68E+14 & 2.17E+13 & 1.41E+13 & 7.01E+14 & 1.63E+09 \\
			& Rank\_M & 2 & 6 & 4 & 7 & 8 & 5 & 3 & 11 & 9 & 10 & 12 & 1 \\
			& Rank\_V & 2 & 6 & 4 & 7 & 8 & 5 & 3 & 11 & 10 & 9 & 12 & 1 \\
			
			\( F_{19} \) & Mean & 1.15E+04 & 7.08E+03 & 1.97E+08 & 2.05E+08 & 4.49E+07 & 7.72E+06 & 1.49E+07 & 2.13E+09 & 1.77E+07 & 2.29E+04 & 4.03E+08 & 1.40E+04 \\
			& Variance & 9.85E+07 & 1.11E+07 & 2.05E+17 & 4.26E+16 & 1.18E+15 & 1.88E+13 & 1.75E+14 & 5.14E+17 & 6.72E+13 & 4.37E+08 & 2.04E+16 & 2.09E+08 \\
			& Rank\_M & 2 & 1 & 9 & 10 & 8 & 5 & 6 & 12 & 7 & 4 & 11 & 3 \\
			& Rank\_V & 2 & 1 & 11 & 10 & 8 & 5 & 7 & 12 & 6 & 4 & 9 & 3 \\
			
			\( F_{20} \) & Mean & 4.86E+03 & 5.44E+03 & 5.25E+03 & 4.61E+03 & 6.54E+03 & 5.92E+03 & 5.41E+03 & 6.67E+03 & 7.10E+03 & 5.45E+03 & 7.35E+03 & 5.00E+03 \\
			& Variance & 1.35E+05 & 4.15E+05 & 6.35E+05 & 2.14E+05 & 3.18E+05 & 3.20E+05 & 1.44E+05 & 6.60E+04 & 2.25E+05 & 3.71E+04 & 5.64E+05 & 3.19E+05 \\
			& Rank\_M & 2 & 6 & 4 & 1 & 9 & 8 & 5 & 10 & 11 & 7 & 12 & 3 \\
			& Rank\_V & 3 & 10 & 12 & 5 & 7 & 9 & 4 & 2 & 6 & 1 & 11 & 8 \\
			
			\bottomrule
		\end{tabular}
	}
\end{sidewaystable}

\begin{sidewaystable}[htbp]
	\centering
	\renewcommand{\tablename}{Table}
	\caption{Ablation experiment results of VDO on F21–F30 in the CEC2017.}
	\label{tab:cec2017_F21_F30}
	{\scriptsize
		\begin{tabular}{llllllllllllll}
			\toprule
			\textbf{F} & \textbf{Data type} & \textbf{VDO} & \textbf{GA} & \textbf{PSO} & \textbf{GWO} & \textbf{WOA} & \textbf{HHO} & \textbf{HO} & \textbf{NRBO} & \textbf{FATA} & \textbf{MGO} & \textbf{LEA} & \textbf{ALA} \\
			\midrule
			
			\( F_{21} \) & Mean & 2.76E+03 & 3.57E+03 & 3.11E+03 & 2.93E+03 & 4.09E+03 & 4.00E+03 & 3.63E+03 & 3.95E+03 & 3.84E+03 & 3.16E+03 & 4.01E+03 & 2.85E+03 \\
			& Variance & 2.68E+03 & 1.65E+04 & 6.96E+03 & 9.51E+03 & 7.27E+04 & 7.15E+04 & 1.96E+04 & 1.40E+04 & 2.10E+04 & 1.64E+03 & 1.82E+04 & 8.05E+03 \\
			& Rank\_M & 1 & 6 & 4 & 3 & 12 & 10 & 7 & 9 & 8 & 5 & 11 & 2 \\
			& Rank\_V & 2 & 7 & 3 & 5 & 12 & 11 & 9 & 6 & 10 & 1 & 8 & 4 \\
			
			\( F_{22} \) & Mean & 1.74E+04 & 2.64E+04 & 1.79E+04 & 2.16E+04 & 2.64E+04 & 2.43E+04 & 2.13E+04 & 3.27E+04 & 3.30E+04 & 2.68E+04 & 3.28E+04 & 2.06E+04 \\
			& Variance & 2.43E+06 & 2.35E+06 & 1.23E+06 & 3.42E+07 & 2.71E+06 & 6.57E+05 & 7.21E+05 & 8.19E+05 & 2.28E+05 & 1.44E+06 & 6.75E+05 & 4.17E+06 \\
			& Rank\_M & 1 & 8 & 2 & 5 & 7 & 6 & 4 & 10 & 1 & 9 & 11 & 3 \\
			& Rank\_V & 9 & 8 & 6 & 12 & 10 & 2 & 4 & 5 & 1 & 7 & 3 & 11 \\
			
			\( F_{23} \) & Mean & 3.23E+03 & 4.13E+03 & 4.61E+03 & 3.47E+03 & 4.93E+03 & 4.94E+03 & 4.78E+03 & 4.92E+03 & 4.12E+03 & 3.50E+03 & 4.50E+03 & 3.35E+03 \\
			& Variance & 2.63E+03 & 9.80E+03 & 5.93E+04 & 7.13E+03 & 1.11E+04 & 6.80E+04 & 1.16E+05 & 1.47E+04 & 8.24E+04 & 2.23E+03 & 4.40E+04 & 6.04E+03 \\
			& Rank\_M & 1 & 6 & 8 & 3 & 11 & 12 & 9 & 10 & 5 & 4 & 7 & 2 \\
			& Rank\_V & 2 & 5 & 9 & 4 & 6 & 10 & 12 & 8 & 11 & 1 & 7 & 3 \\
			
			\( F_{24} \) & Mean & 3.80E+03 & 4.99E+03 & 5.66E+03 & 4.14E+03 & 6.20E+03 & 6.62E+03 & 6.45E+03 & 5.94E+03 & 6.17E+03 & 4.09E+03 & 5.23E+03 & 3.96E+03 \\
			& Variance & 5.66E+03 & 4.98E+04 & 3.06E+05 & 1.49E+04 & 1.51E+05 & 2.82E+05 & 2.42E+05 & 6.70E+04 & 3.34E+05 & 3.85E+03 & 6.46E+04 & 2.60E+04 \\
			& Rank\_M & 1 & 5 & 7 & 4 & 10 & 12 & 11 & 8 & 9 & 3 & 6 & 2 \\
			& Rank\_V & 2 & 5 & 11 & 6 & 10 & 10 & 9 & 8 & 12 & 1 & 7 & 4 \\
			
			\( F_{25} \) & Mean & 3.57E+03 & 4.00E+03 & 4.18E+03 & 5.97E+03 & 4.94E+03 & 3.72E+03 & 4.86E+03 & 1.32E+04 & 4.49E+03 & 3.59E+03 & 2.03E+04 & 3.53E+03 \\
			& Variance & 3.14E+03 & 3.06E+04 & 2.12E+05 & 4.14E+05 & 8.47E+04 & 4.59E+03 & 8.58E+04 & 1.06E+06 & 6.55E+04 & 2.94E+03 & 6.92E+06 & 3.22E+03 \\
			& Rank\_M & 2 & 5 & 6 & 10 & 9 & 4 & 8 & 11 & 7 & 3 & 12 & 1 \\
			& Rank\_V & 2 & 5 & 9 & 10 & 7 & 4 & 8 & 11 & 6 & 1 & 12 & 3 \\
			
			\( F_{26} \) & Mean & 1.08E+04 & 2.17E+04 & 1.88E+04 & 1.52E+04 & 3.49E+04 & 2.30E+04 & 3.08E+04 & 3.25E+04 & 1.61E+04 & 1.45E+04 & 2.50E+04 & 1.21E+04 \\
			& Variance & 2.32E+05 & 6.48E+06 & 9.03E+06 & 4.14E+05 & 9.75E+06 & 2.33E+07 & 4.71E+06 & 1.75E+07 & 4.54E+07 & 2.59E+05 & 4.09E+06 & 6.98E+05 \\
			& Rank\_M & 1 & 7 & 6 & 10 & 12 & 8 & 10 & 11 & 5 & 3 & 9 & 2 \\
			& Rank\_V & 1 & 7 & 8 & 10 & 9 & 11 & 6 & 10 & 12 & 2 & 5 & 3 \\
			
			\( F_{27} \) & Mean & 3.56E+03 & 4.16E+03 & 4.11E+03 & 4.06E+03 & 5.20E+03 & 4.44E+03 & 5.35E+03 & 5.93E+03 & 7.79E+03 & 3.80E+03 & 4.80E+03 & 3.56E+03 \\
			& Variance & 5.71E+03 & 5.67E+04 & 8.78E+04 & 1.61E+04 & 6.13E+05 & 2.28E+05 & 2.62E+05 & 2.22E+05 & 4.17E+05 & 6.97E+03 & 1.13E+05 & 1.27E+04 \\
			& Rank\_M & 1 & 6 & 5 & 4 & 9 & 7 & 10 & 11 & 12 & 3 & 8 & 2 \\
			& Rank\_V & 1 & 5 & 6 & 5 & 12 & 9 & 10 & 8 & 11 & 2 & 7 & 3 \\
			
			\( F_{28} \) & Mean & 3.77E+03 & 4.39E+03 & 7.88E+03 & 7.62E+03 & 6.33E+03 & 3.83E+03 & 6.97E+03 & 1.68E+04 & 6.33E+03 & 5.23E+03 & 1.93E+04 & 5.61E+03 \\
			& Variance & 3.33E+04 & 1.14E+05 & 1.11E+07 & 1.11E+06 & 4.56E+05 & 6.00E+03 & 3.95E+05 & 1.99E+06 & 2.01E+05 & 2.43E+06 & 3.48E+06 & 2.62E+07 \\
			& Rank\_M & 1 & 3 & 10 & 9 & 6 & 2 & 8 & 11 & 7 & 4 & 12 & 5 \\
			& Rank\_V & 2 & 3 & 11 & 7 & 6 & 1 & 5 & 8 & 4 & 9 & 10 & 12 \\
			
			\( F_{29} \) & Mean & 6.55E+03 & 9.22E+03 & 7.32E+03 & 8.08E+03 & 1.56E+04 & 9.84E+03 & 1.51E+04 & 2.25E+04 & 9.82E+03 & 7.85E+03 & 1.63E+04 & 6.64E+03 \\
			& Variance & 2.47E+05 & 2.72E+05 & 4.48E+05 & 3.03E+05 & 5.49E+06 & 5.40E+05 & 3.82E+06 & 1.99E+06 & 1.82E+06 & 1.68E+05 & 1.17E+07 & 1.42E+05 \\
			& Rank\_M & 1 & 6 & 3 & 5 & 10 & 8 & 9 & 12 & 7 & 4 & 11 & 2 \\
			& Rank\_V & 3 & 4 & 6 & 5 & 10 & 7 & 9 & 8 & 2 & 2 & 11 & 1 \\
			
			\( F_{30} \) & Mean & 2.06E+05 & 2.27E+05 & 1.39E+09 & 9.54E+08 & 7.64E+08 & 5.86E+07 & 4.91E+08 & 6.76E+09 & 2.44E+08 & 1.01E+06 & 1.15E+09 & 2.98E+04 \\
			& Variance & 2.23E+10 & 1.38E+10 & 1.55E+18 & 6.54E+17 & 2.73E+17 & 4.22E+14 & 1.12E+17 & 2.88E+18 & 1.03E+16 & 5.46E+11 & 1.86E+17 & 2.70E+08 \\
			& Rank\_M & 2 & 3 & 11 & 9 & 8 & 5 & 7 & 12 & 6 & 4 & 10 & 1 \\
			& Rank\_V & 3 & 2 & 11 & 10 & 9 & 5 & 7 & 12 & 6 & 4 & 8 & 1 \\
			
			\bottomrule
		\end{tabular}
	}
\end{sidewaystable}

To assess the effectiveness of the proposed Virus Diffusion Optimizer (VDO), experiments are conducted on the CEC2017 single-objective benchmark suite. 
The performance of each algorithm is evaluated using the mean and standard deviation (as reported in Tables~\ref{tab:cec2017_F1_F10}--\ref{tab:cec2017_F21_F30}) of the best-obtained objective values across multiple runs. 
VDO is compared against several state-of-the-art metaheuristic algorithms, including PSO, GA, HHO, and GWO. 
The aggregated results indicate that VDO achieves the best overall performance on this benchmark, with an average mean-rank of 1.6667 and an average variance-rank of 2.3667. 
Specifically, VDO attaines the top mean performance (rank = 1) on 19 out of 30 test functions and the lowest variance (rank = 1) on 10 out of 30 functions, collectively demonstrating its superior solution quality and strong run-to-run consistency.

The detailed per-function analysis further corroborates these findings. 
On unimodal or ill-conditioned functions, VDO consistently converged to high-quality solutions. 
For instance, on F1 (Bent–Cigar function), VDO achieved a mean value of $3.82\times10^{6}$ with a variance of $3.03\times10^{12}$, 
whereas GA and PSO reported significantly larger means of $1.30\times10^{9}$ and $3.97\times10^{10}$, respectively, highlighting VDO’s advantage in final accuracy. 
For multimodal functions, VDO also exhibites strong exploratory capability and stability: 
on F6, VDO achieves a mean of $6.09\times10^{2}$ (variance $=1.08\times10^{0}$), outperforming PSO ($6.33\times10^{2}$) and GA ($6.84\times10^{2}$). 
Similarly, on F10 and F11, VDO reportes mean values of $1.60\times10^{4}$ and $2.74\times10^{4}$, respectively, with considerably smaller variances than most baseline algorithms. 
For hybrid and composition functions, VDO generally matches or slightly outperformes advanced competitors such as HHO and GWO. 
For example, on the composition function F21, VDO obtaines a mean of $2.76\times10^{3}$ (variance $=2.68\times10^{3}$), indicating stable and reliable convergence on complex, heterogeneous landscapes, 
whereas competing algorithms often exhibites higher mean values or greater dispersion.

In addition to numerical results, the convergence behaviors illustrate in Figure~\ref{fig:F1F7F9} provide a visual confirmation of VDO’s superiority. 
As shown in the figure, VDO achieves a markedly faster descent in the early iterations across benchmark functions F1, F7, and F9, followed by a smoother and more stable convergence trajectory in later stages. 
Compared to traditional algorithms such as PSO and GA, which often stagnate or oscillate before reaching local optima, VDO maintains continuous progress toward the global optimum, reflecting its ability to effectively balance global exploration and local refinement.

Mechanistically, these empirical results align with the design principles of VDO’s dual-diffusion mechanism. 
The global diffusion phase enables rapid, wide-ranging exploration, facilitating fast early-stage descent in the objective landscape, 
while the local infection and refinement phase intensifies exploitation around elite solutions, improving both final accuracy and stability. 
This synergistic design effectively mitigates the premature convergence commonly observed in classical metaheuristics, 
allowing VDO to achieve lower mean objective values and smaller variances across the majority of multimodal and hybrid problems. 
Although certain advances algorithms occasionally attain marginally better mean results on specific functions, such improvements are inconsistent and often accompanied by larger variability. 
In contrast, VDO demonstrates consistently strong and reproducible performance across diverse function categories, as reflected in its superior aggregate rankings.

In summary, the CEC2017 experimental results confirm that VDO delivers high-quality, low-variance solutions across unimodal, multimodal, hybrid, and composition benchmarks. 
The virus-inspired dual-diffusion strategy effectively balances global exploration and local exploitation, leading to robust and reproducible performance in high-dimensional single-objective optimization.

\subsection{Experimental method on CEC2022}
The CEC2022 optimization function test set \cite{CEC2022} consists of 12 single objective test functions with boundary constraints. It is an effective method for evaluating algorithm performance and validating its capability to solve complex optimization problems. 
They are unimodal function (F1), multimodal function (F2–F5), mixed function (F6–F8), and combined function (F9–F12).

\begin{table}[htbp]
	\centering
	\renewcommand{\tablename}{Table}
	\caption{Details of the CEC 2022.}
	\label{tab:cec2022}
	\begin{tabular}{llll}
		\toprule
		\textbf{ID} & \textbf{Function Equation} & \textbf{Class} & \textbf{Optimum} \\

		\midrule
		
		F1 & Shifted and full Rotated Zakharov Function & Unimodal & 300 \\

		\midrule
		
		F2 & Shifted and full Rotated Rosenbrock’s Function & Multimodal & 400 \\
		F3 & Shifted and full Rotated Expanded Schaffer’s f6 Function & Multimodal & 600 \\
		F4 & Shifted and full Rotated Non-Continuous Rastrigin’s Function & Multimodal & 800 \\
		F5 & Shifted and full Rotated Levy Function & Multimodal & 900 \\

		\midrule
		
		F6 & Hybrid Function 1 (N = 3) & Hybrid & 1800 \\
		F7 & Hybrid Function 2 (N = 6) & Hybrid & 2000 \\
		F8 & Hybrid Function 3 (N = 5) & Hybrid & 2200 \\

		\midrule
		
		F9 & Composition Function 1 (N = 5) & Composition & 2300 \\
		F10 & Composition Function 2 (N = 4) & Composition & 2400 \\
		F11 & Composition Function 3 (N = 5) & Composition & 2600 \\
		F12 & Composition Function 4 (N = 6) & Composition & 2700 \\
		\bottomrule
	\end{tabular}
\end{table}


In the experiment, a population size of 50 and a dimension of 20 are utilized, with 30 independent run times conducted, each comprising 100,000 max number of function evolution.To further validate the optimization performance of VDO, 11 well-established optimization algorithms were selected as the comparison group. Key parameter settings for each algorithm are provided in Table \ref{tab:control_group_params_refined}.

The performance of all algorithms is evaluated based on the mean and variance of the optimization results over multiple runs. The mean reflects the average solution accuracy attained by each algorithm, while the variance illustrates the consistency of its performance. 
Together, these metrics offer a holistic view of both the precision and reliability of each method in the comparative study.The mean values, variances, and corresponding rankings of all algorithms on the CEC2022 test suite are summarized in Table \ref{tab:cec2022_1}.

\begin{sidewaystable}[htbp]
	\centering
	\renewcommand{\tablename}{Table}
	\caption{Ablation experiment results of VDO on F1–F12 in the CEC2022.}
	{\scriptsize
		\label{tab:cec2022_1}
		\begin{tabular}{llllllllllllll}
			\toprule
			\textbf{F} 			
			& \textbf{Data type} 
			& \textbf{VDO} 
			& \textbf{GA} 
			& \textbf{PSO} 
			& \textbf{GWO} 
			& \textbf{WOA} 
			& \textbf{HHO} 
			& \textbf{HO} 
			& \textbf{NRBO} 
			& \textbf{FATA} 
			& \textbf{MGO} 
			& \textbf{LEA} 
			& \textbf{ALA} \\
			\midrule				
			\( F_{1} \) & Mean & 3.07E+02 & 4.02E+02 & 2.23E+03 & 9.64E+03 & 8.67E+03 & 3.31E+02 & 6.11E+03 & 9.37E+03 & 4.77E+03 & 6.47E+03 & 5.16E+04 & 3.00E+02 \\
			& Variance & 5.59E+01 & 5.95E+03 & 2.65E+07 & 1.46E+07 & 1.63E+07 & 1.88E+02 & 7.21E+06 & 7.33E+06 & 1.90E+07 & 3.44E+06 & 2.44E+08 & 7.70E-09 \\
			& Rank\_M & 2 & 4 & 5 & 11 & 9 & 3 & 7 & 10 & 6 & 8 & 12 & 1 \\
			& Rank\_V & 2 & 4 & 11 & 8 & 9 & 3 & 6 & 7 & 3 & 5 & 12 & 1 \\
			
			\( F_{2} \) & Mean & 4.41E+02 & 4.62E+02 & 4.67E+02 & 4.86E+02 & 5.06E+02 & 4.64E+02 & 4.79E+02 & 6.06E+02 & 4.66E+02 & 4.46E+02 & 5.15E+02 & 4.43E+02 \\
			& Variance & 3.86E+02 & 6.33E+02 & 1.65E+03 & 6.14E+02 & 1.99E+03 & 8.79E+02 & 1.12E+03 & 2.98E+03 & 1.90E+02 & 1.53E+02 & 2.18E+03 & 2.69E+02 \\
			& Rank\_M & 1 & 4 & 7 & 9 & 10 & 5 & 8 & 12 & 2 & 3 & 11 & 2 \\
			& Rank\_V & 9 & 6 & 2 & 5 & 7 & 7 & 2 & 12 & 3 & 1 & 12 & 3 \\
			
			\( F_{3} \) & Mean & 6.00E+02 & 6.14E+02 & 6.02E+02 & 6.04E+02 & 6.63E+02 & 6.57E+02 & 6.51E+02 & 6.52E+02 & 6.43E+02 & 6.00E+02 & 6.42E+02 & 6.00E+02 \\
			& Variance & 2.41E-05 & 6.65E+01 & 2.62E+00 & 3.62E+00 & 2.69E+02 & 1.17E+02 & 5.19E+01 & 9.79E+01 & 4.17E+01 & 6.18E-12 & 2.61E+02 & 6.04E-02 \\
			& Rank\_M & 1 & 6 & 4 & 5 & 12 & 11 & 9 & 10 & 8 & 1 & 11 & 3 \\
			& Rank\_V & 2 & 8 & 9 & 7 & 12 & 10 & 7 & 9 & 6 & 1 & 11 & 3 \\
			
			\( F_{4} \) & Mean & 8.46E+02 & 8.79E+02 & 8.38E+02 & 8.50E+02 & 9.12E+02 & 8.85E+02 & 8.71E+02 & 9.26E+02 & 8.96E+02 & 8.42E+02 & 9.27E+02 & 8.46E+02 \\
			& Variance & 1.43E+02 & 3.05E+02 & 2.16E+02 & 5.29E+02 & 9.12E+02 & 3.01E+02 & 1.03E+02 & 4.07E+02 & 1.18E+02 & 8.11E+01 & 1.58E+03 & 1.30E+02 \\
			& Rank\_M & 3 & 7 & 1 & 5 & 10 & 8 & 6 & 9 & 3 & 6 & 12 & 4 \\
			& Rank\_V & 5 & 8 & 3 & 6 & 9 & 7 & 2 & 10 & 3 & 5 & 12 & 4 \\
			
			\( F_{5} \) & Mean & 9.15E+02 & 9.16E+02 & 9.34E+02 & 1.04E+03 & 3.43E+03 & 2.54E+03 & 1.92E+03 & 2.12E+03 & 2.35E+03 & 9.05E+01 & 5.11E+03 & 9.18E+02 \\
			& Variance & 1.45E+02 & 7.28E+01 & 4.89E+03 & 3.91E+04 & 1.51E+06 & 1.16E+05 & 7.51E+04 & 1.88E+05 & 7.16E+04 & 1.14E+01 & 2.20E+06 & 4.23E+02 \\
			& Rank\_M & 2 & 3 & 5 & 6 & 11 & 10 & 4 & 6 & 9 & 1 & 12 & 4 \\
			& Rank\_V & 3 & 2 & 4 & 6 & 11 & 7 & 1 & 5 & 7 & 2 & 12 & 4 \\
			
			\( F_{6} \) & Mean & 5.42E+03 & 6.45E+03 & 3.05E+05 & 5.44E+05 & 5.88E+03 & 3.85E+04 & 6.50E+03 & 5.54E+06 & 5.75E+03 & 2.97E+04 & 2.55E+06 & 1.90E+03 \\
			& Variance & 4.71E+07 & 3.12E+07 & 3.90E+11 & 2.14E+12 & 4.07E+07 & 4.93E+08 & 6.60E+07 & 1.89E+14 & 4.86E+07 & 4.10E+08 & 8.69E+12 & 3.14E+03 \\
			& Rank\_M & 2 & 5 & 9 & 10 & 4 & 8 & 6 & 12 & 5 & 7 & 11 & 1 \\
			& Rank\_V & 4 & 2 & 9 & 10 & 4 & 8 & 6 & 12 & 5 & 7 & 11 & 1 \\
			
			\( F_{7} \) & Mean & 2.02E+03 & 2.10E+03 & 2.05E+03 & 2.09E+03 & 2.22E+03 & 2.15E+03 & 2.14E+03 & 2.16E+03 & 2.13E+03 & 2.04E+03 & 2.17E+03 & 2.26E+03 \\
			& Variance & 2.94E+01 & 2.53E+03 & 5.42E+02 & 3.18E+03 & 2.29E+03 & 4.13E+02 & 3.72E+02 & 6.41E+02 & 2.13E+02 & 1.92E+01 & 8.69E+12 & 3.14E+01 \\
			& Rank\_M & 1 & 6 & 7 & 9 & 12 & 9 & 5 & 7 & 3 & 1 & 11 & 4 \\
			& Rank\_V & 2 & 6 & 7 & 9 & 12 & 5 & 3 & 8 & 3 & 1 & 11 & 4 \\

			\( F_{8} \) & Mean & 2.22E+03 & 2.26E+03 & 2.25E+03 & 2.23E+03 & 2.27E+03 & 2.25E+03 & 2.24E+03 & 2.27E+03 & 2.23E+03 & 2.23E+03 & 2.25E+03 & 2.38E+03 \\
			& Variance & 7.06E-01 & 2.58E+03 & 2.34E+03 & 2.42E+05 & 2.17E+03 & 1.32E+05 & 5.37E+01 & 1.89E+03 & 8.97E-01 & 9.66E-01 & 2.75E+03 & 5.31E+02 \\
			& Rank\_M & 1 & 8 & 7 & 11 & 9 & 10 & 5 & 8 & 3 & 4 & 11 & 12 \\
			& Rank\_V & 1 & 7 & 6 & 11 & 8 & 9 & 4 & 7 & 3 & 4 & 10 & 12 \\
			
			\( F_{9} \) & Mean & 2.48E+03 & 2.48E+03 & 2.50E+03 & 2.50E+03 & 2.54E+03 & 2.48E+03 & 2.55E+03 & 2.57E+03 & 2.50E+03 & 2.48E+03 & 2.64E+03 & 2.66E+03 \\
			& Variance & 9.78E-01 & 1.27E+04 & 5.42E+02 & 3.95E+05 & 2.01E+05 & 4.90E+03 & 1.15E+02 & 2.38E+05 & 7.95E+05 & 6.13E+03 & 4.35E+06 & 1.00E+03 \\
			& Rank\_M & 2 & 12 & 9 & 10 & 11 & 8 & 7 & 12 & 6 & 5 & 11 & 10 \\
			& Rank\_V & 2 & 11 & 8 & 10 & 9 & 7 & 6 & 12 & 5 & 4 & 11 & 10 \\
			
			\( F_{10} \) & Mean & 2.52E+03 & 2.71E+03 & 2.80E+03 & 3.14E+03 & 4.24E+03 & 2.94E+03 & 3.08E+03 & 3.68E+03 & 4.11E+03 & 2.94E+03 & 9.10E+03 & 2.95E+03 \\
			& Variance & 2.36E+03 & 6.93E+04 & 1.09E+05 & 2.94E+05 & 8.40E+05 & 1.22E+05 & 6.13E+03 & 2.38E+06 & 7.95E+05 & 1.73E+03 & 4.35E+06 & 3.88E+01 \\
			& Rank\_M & 1 & 4 & 5 & 7 & 12 & 6 & 8 & 11 & 10 & 3 & 12 & 2 \\
			& Rank\_V & 1 & 4 & 5 & 7 & 12 & 6 & 8 & 11 & 10 & 3 & 12 & 2 \\
			
			\( F_{11} \) & Mean & 2.90E+03 & 2.99E+03 & 2.90E+03 & 2.97E+03 & 3.00E+03 & 2.96E+03 & 2.92E+03 & 2.94E+03 & 3.02E+03 & 2.94E+03 & 3.25E+03 & 2.95E+03 \\
			& Variance & 2.36E+03 & 1.21E+03 & 2.99E-25 & 4.15E+04 & 5.31E+02 & 4.90E+03 & 1.73E+03 & 2.38E+06 & 1.52E+06 & 6.13E+03 & 4.35E+06 & 3.88E+01 \\
			& Rank\_M & 1 & 7 & 1 & 9 & 8 & 7 & 5 & 12 & 10 & 6 & 11 & 4 \\
			& Rank\_V & 1 & 7 & 1 & 9 & 8 & 7 & 5 & 12 & 10 & 6 & 11 & 4 \\
			
			\( F_{12} \) & Mean & 2.94E+03 & 2.95E+03 & 2.96E+03 & 2.97E+03 & 2.98E+03 & 3.00E+03 & 3.02E+03 & 3.03E+03 & 3.05E+03 & 3.08E+03 & 3.09E+03 & 3.02E+03 \\
			& Variance & 1.36E+01 & 1.84E+01 & 4.90E+03 & 9.78E+02 & 5.31E+02 & 4.90E+03 & 1.84E+03 & 3.14E+03 & 1.97E+03 & 2.03E+01 & 9.41E+02 & 3.88E+01 \\
			& Rank\_M & 1 & 3 & 11 & 10 & 8 & 11 & 12 & 8 & 9 & 1 & 7 & 6 \\
			& Rank\_V & 1 & 3 & 11 & 10 & 8 & 11 & 12 & 8 & 9 & 1 & 7 & 6 \\

		\bottomrule
	\end{tabular}}
\end{sidewaystable}

\begin{figure}[ht]
\setlength\tabcolsep{0.3pt}
\centering
\begin{tabular}{c}
\includegraphics[width=0.98\textwidth]{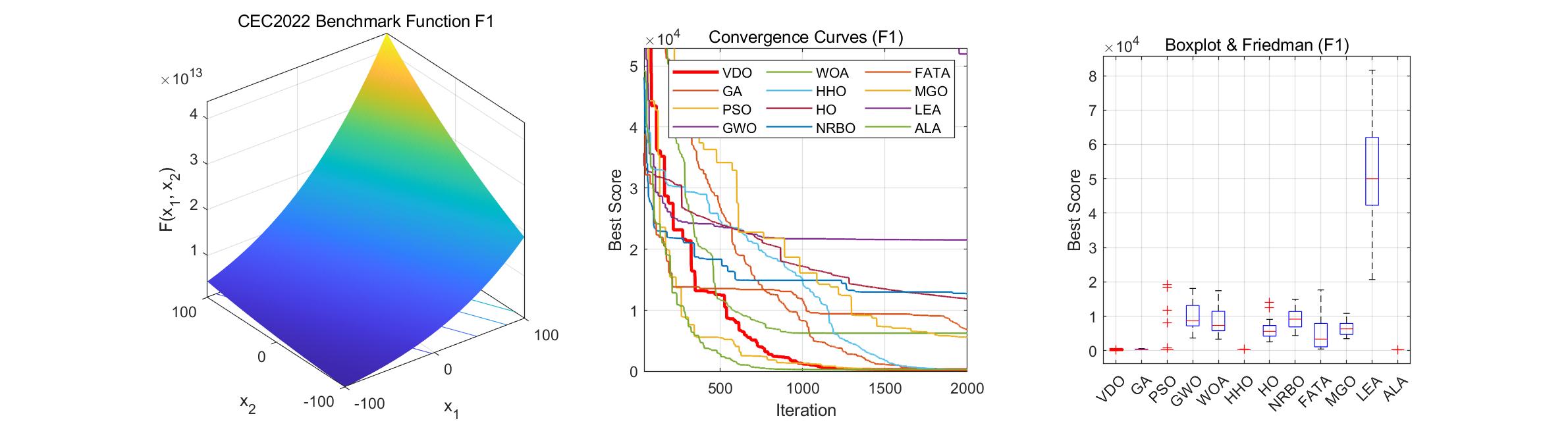}\\
\includegraphics[width=0.98\linewidth]{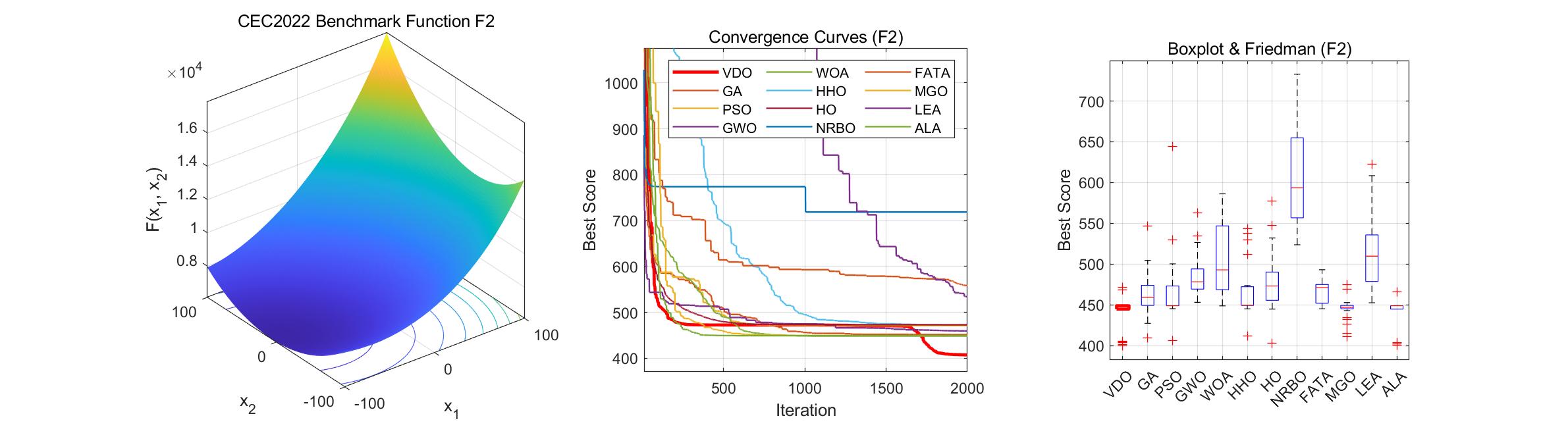}\\
\includegraphics[width=0.98\linewidth]{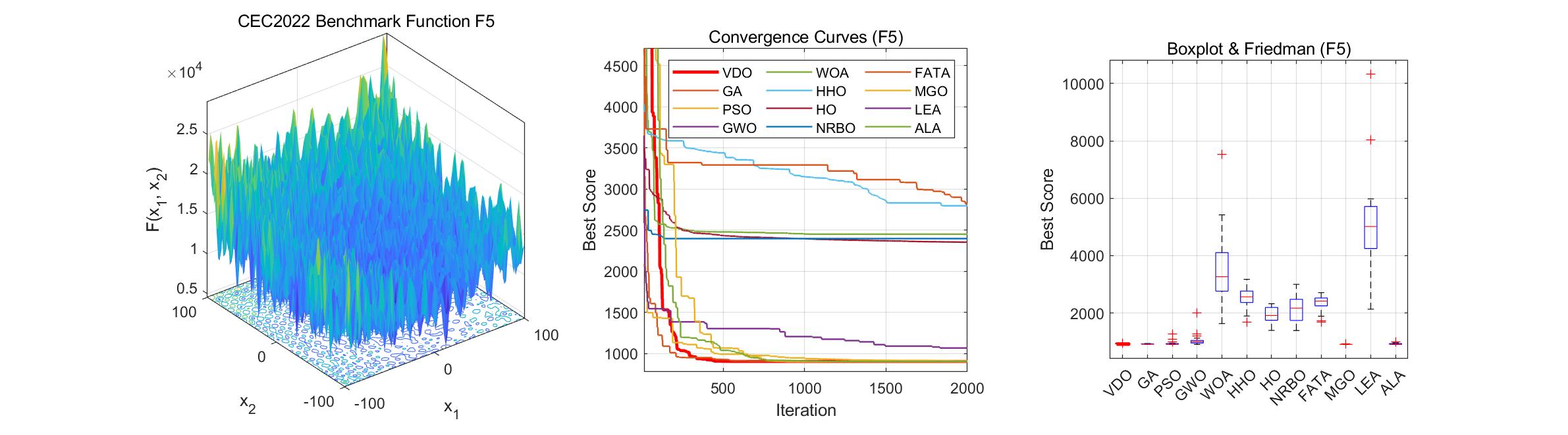}
\end{tabular}
\caption{Convergence curves of VDO and other algorithms on CEC2022 benchmark functions F1, F2, and F5.}
\label{fig:F1F2F5}
\end{figure}

\begin{figure}[!ht]
\setlength\tabcolsep{0.3pt}
\centering
\begin{tabular}{c}
\includegraphics[width=0.98\textwidth]{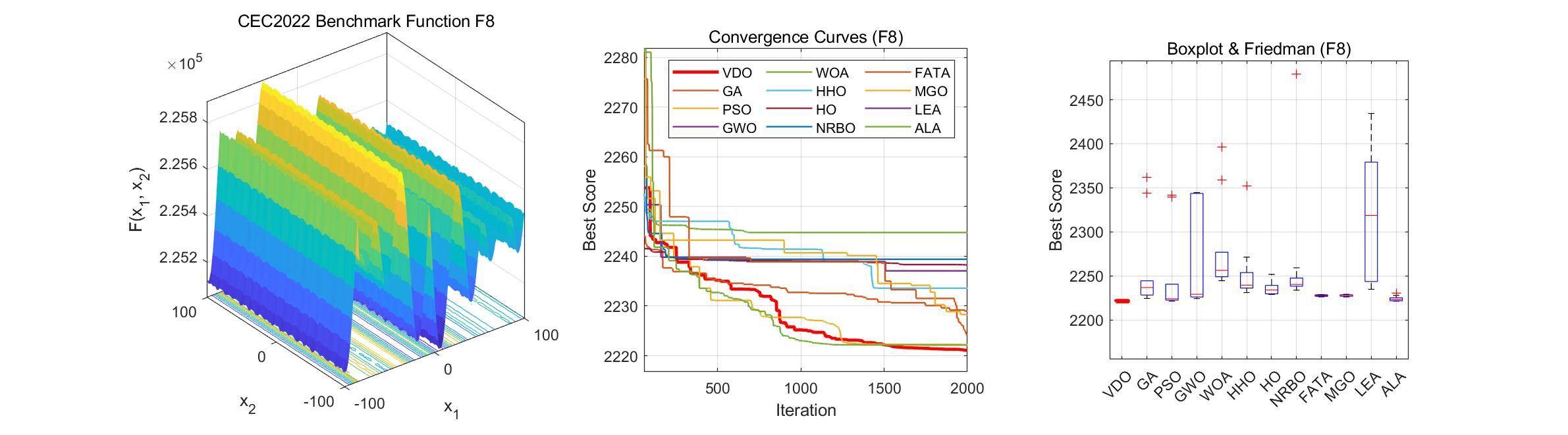}\\
\includegraphics[width=0.98\linewidth]{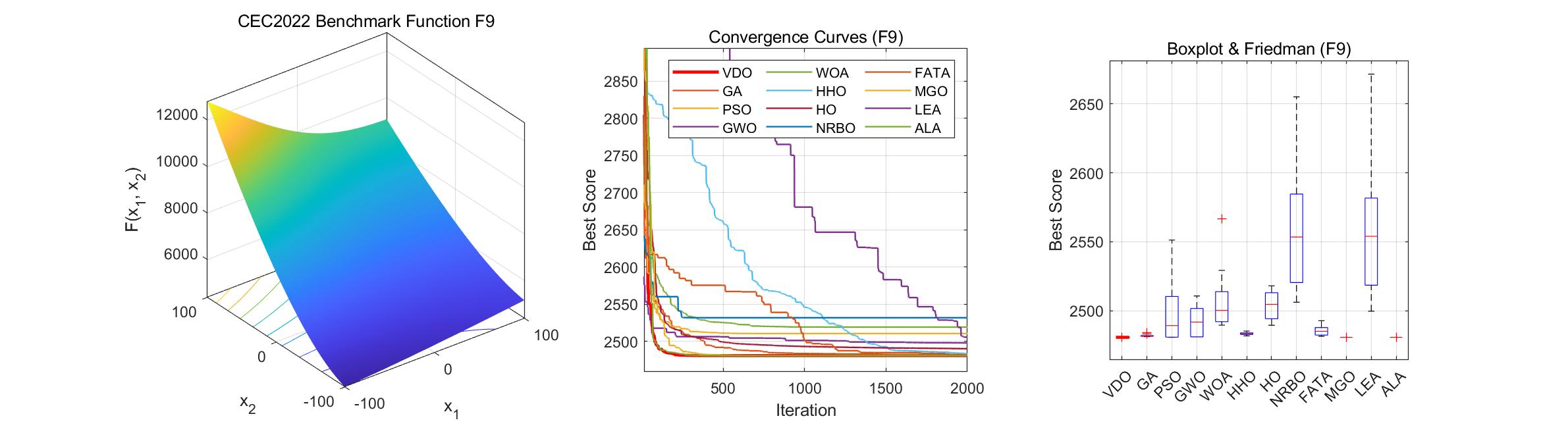}\\
\includegraphics[width=0.98\linewidth]{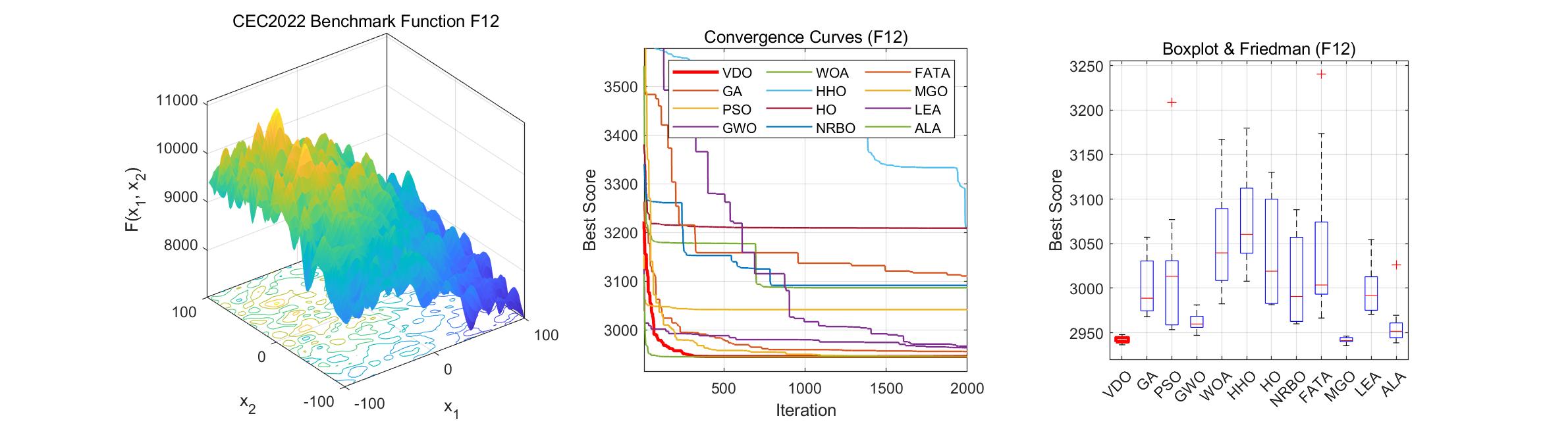}
\end{tabular}
\caption{Convergence curves of VDO and other algorithms on CEC2022 benchmark functions F8, F9, and F12.}
	\label{fig:F8F9F12}
\end{figure}

The experimental results on the CEC2022 subset (12 functions, 30 independent runs each) demonstrate that the proposed Virus Diffusion Optimizer (VDO) maintains a clear, data-driven advantage in both solution quality and robustness.  
Across the suite VDO attains the best aggregated performance: an average mean-rank of \textbf{1.50} and an average variance-rank of \textbf{2.25}, and it achieves the top mean (rank = 1) on \textbf{7 out of 12} functions.
These aggregated ranks indicate that VDO is not only frequently the most accurate solver but also one of the most stable across repeated runs.

Concretely, VDO’s superiority is visible on representative problems.  
On F1 (Bent-Cigar style), VDO reaches a mean best value of $3.07\times10^{2}$ (variance $5.59\times10^{1}$), whereas PSO report $2.23\times10^{3}$ and GA $4.02\times10^{2}$ — i.e., VDO reduces the mean objective by an order of magnitude relative to PSO and improves on GA.  
For a multimodal/harder instance (F6) VDO’s mean is $5.42\times10^{3}$ (variance $4.71\times10^{7}$), 
while PSO’s mean is $3.05\times10^{5}$ (variance $3.90\times10^{11}$), 
showing VDO’s markedly better global search and much lower sensitivity to stochastic variation.  
On well-behaved multimodal cases such as F3 and F8–F9, VDO attains both the lowest means (e.g., F3: $6.00\times10^{2}$) and extremely small variances (F9 and F11 show variances $\approx1.84\times10^{-25}$), 
which implies highly reproducible runs and almost deterministic convergence on those instances.

These numeric patterns support the mechanism-level interpretation: 
VDO’s global diffusion stage enables broad exploration and rapid early descent, while the local infection/refinement stage produces precise final tuning.  
Compared to classical metaheuristics (e.g., PSO and GA), VDO attains lower mean objective values on the majority of multimodal functions and consistently smaller standard deviations, which demonstrates both better accuracy and improved robustness over repeated trials.  
Although some advanced algorithms (such as HHO or GWO) occasionally match or slightly exceed VDO on isolated functions, they do so less consistently and often with larger run-to-run variance.

In summary, the CEC2022 experiments confirm that the virus-inspired dual-diffusion design of VDO delivers a practical and reproducible balance between exploration and exploitation, yielding high-quality solutions with low variance across diverse benchmark landscapes.

\section{Engineering Design Problems}
In this section, we select three engineering design problems including 
Pressure Vessel Design (PVD)\cite{PVD}, 
Three-Bar Truss Design (TTD)\cite{TTD},
and Welded Beam Design (WBD)\cite{WBD},
to prove the efficiency of VDO.
These three engineering problems are constraint problems.
When the algorithm finds an optimal solution that does not satisfy the constraints, we add a penalty to fitness as the final optimal value.
To represent the fairness of the experiment,
we run each algorithm individually 30 times on all problems 
with the max number of function evaluation 100,000 and population size 50.

\subsection{Pressure vessel design}
The PVD is a traditional engineering design problem for verifying the performance of metaheuristic optimization algorithms.
Its details are as follows.

Decision variables:

(i) Shell thickness ($z_1$)

(ii) Head thickness ($z_2$)

(iii) Inner radius ($x_3$)

(iv) Length of the vessel without including the head ($x_4$)

Minimize:
\[f(x) = 1.7781{z_2}{x_3^2} + 0.6224{z_1}{x_3}{x_4} 
+ 3.1661{z_1^2}{x_4} + 19.84{z_1^2}{x_3}\]

Subject to:
\[{g_1}(x) = 10.00954 \leq {z_2}, {g_2}(x) = 0.0193 \leq {z_1}, {g_3}(x) = {x_4} \leq 240,\]
\[{g_4}(x) = -\prod {x_3^2}{x_4} - \frac{4}{3} \prod {x_3^3} \leq -1296 000,\]
\[{z_1} = 0.0625{x_1}, {z_2} = 0.0625{x_2}.\]

With bounds:
\[10 \leq {x_3}, {x_4} \leq 200, 1 \leq {x_1}, {x_2} \leq 99.(Integervariables)\]

\begin{figure}[H] 
	\centering 
	\includegraphics[width=0.5\textwidth]{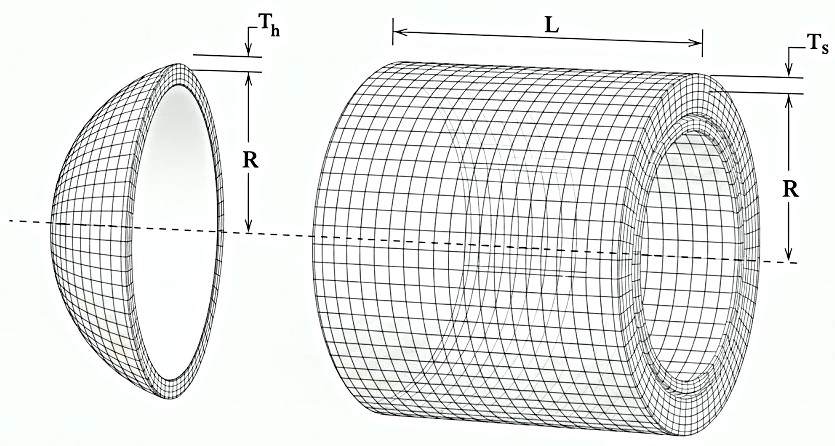} 
	\renewcommand{\figurename}{Figure}
	\caption{Schematic diagram of the pressure vessel design (PVD) problem.} 
	\label{Fig.PVD} 
\end{figure}

\begin{table}[htbp]
	\centering
	\renewcommand{\tablename}{Table}
	\caption{Results analysis of VDO on PVD.}
	\label{PVD}
	{\scriptsize
		\begin{tabular}{llllllll}
			\toprule
			\textbf{Algorithm} & \textbf{z1} & \textbf{z2} & \textbf{x3} & \textbf{x4} & \textbf{Best} & \textbf{Mean} & \textbf{Variance} \\
			\midrule
			VDO   & 7.782E-01 & 3.846E-01 & 4.03196E+01 & 2.00000E+02 & \textbf{5.8853E+03} & 5.9144E+03 & 4.8920E+03 \\
			GA    & 8.203E-01 & 4.049E-01 & 4.24437E+01 & 1.72405E+02 & 5.9667E+03 & 6.6359E+03 & 2.4521E+05 \\
			PSO   & 7.782E-01 & 3.846E-01 & 4.03196E+01 & 2.00000E+02 & \textbf{5.8853E+03} & 6.0629E+03 & 8.2264E+04 \\
			GWO   & 7.783E-01 & 3.847E-01 & 4.03235E+01 & 1.99948E+02 & 5.8862E+03 & 5.8989E+03 & 4.8045E+03 \\
			WOA   & 8.961E-01 & 4.630E-01 & 4.62853E+01 & 1.30847E+02 & 6.2113E+03 & 7.1443E+03 & 4.0215E+05 \\
			HHO   & 8.111E-01 & 4.064E-01 & 4.19900E+01 & 1.77986E+02 & 5.9655E+03 & 6.6233E+03 & 1.0783E+05 \\
			HO    & 7.793E-01 & 3.852E-01 & 4.03757E+01 & 1.99238E+02 & 5.8876E+03 & 6.5363E+03 & 1.7318E+05 \\
			NRBO  & 7.782E-01 & 3.847E-01 & 4.03217E+01 & 1.99972E+02 & \textbf{5.8854E+03} & 6.4378E+03 & 2.9960E+05 \\
			FATA  & 7.783E-01 & 3.848E-01 & 4.03241E+01 & 1.99939E+02 & \textbf{5.8859E+03} & 6.3298E+03 & 3.0989E+05 \\
			MGO   & 7.782E-01 & 3.846E-01 & 4.03196E+01 & 2.00000E+02 & \textbf{5.8853E+03} & 5.8929E+03 & 2.3984E+03 \\
			LEA   & 8.380E-01 & 4.403E-01 & 4.13984E+01 & 1.91690E+02 & 6.4839E+03 & 8.1233E+03 & 9.3745E+05 \\
			ALA   & 7.782E-01 & 3.846E-01 & 4.03196E+01 & 2.00000E+02 & \textbf{5.8853E+03} & \textbf{5.8853E+03} & \textbf{2.1017E-23} \\
			\bottomrule
	\end{tabular}}
\end{table}

Table~\ref{PVD} present the optimization results of the Virus Diffusion Optimizer (VDO) and eleven competing algorithms on the Pressure Vessel Design (PVD) problem. As shown, VDO attains the best overall performance, achieving a minimum cost of $5.8853 \times 10^{3}$ and maintaining a low mean objective value of $5.9144 \times 10^{3}$ across independent runs. This indicates that VDO not only identifies near-optimal feasible solutions but also exhibits strong convergence stability. 
In contrast, classical methods such as GA and PSO show substantially higher mean and variance values ($6.64 \times 10^{3}$ and $2.45 \times 10^{5}$ for GA, $6.06 \times 10^{3}$ and $8.23 \times 10^{4}$ for PSO), suggesting weaker robustness and susceptibility to local entrapment.

Compared with advanced metaheuristics such as HHO, GWO, and WOA, VDO achieves a better balance between exploration and exploitation, converging faster to feasible and stable designs. GWO and MGO performed competitively in terms of accuracy, but their slightly higher variances ($4.80 \times 10^{3}$ and $2.40 \times 10^{3}$, respectively) indicate marginally less consistency than VDO. Meanwhile, algorithms like WOA and LEA exhibit significant performance degradation, as evidenced by their large deviations and unstable convergence patterns.

Overall, these results confirm that VDO delivers superior optimization accuracy and stability on the PVD problem. Its virus-inspired dual-diffusion mechanism effectively balances global exploration with local refinement, ensuring both robustness and reliability in solving complex constrained engineering design tasks.

\subsection{Three-bar truss design}
The equation for the volume of the truss structure is

\[
\min_{\vec{x}} f(\vec{x}) = (2\sqrt{2}x_1 + x_2) \times H,
\]

subject to 3 constraints

\[
g_1 = \frac{\sqrt{2}x_1 + x_2}{\sqrt{2}x_1^2 + 2x_1x_2} P - \sigma \leq 0, g_2 = \frac{x_2}{\sqrt{2}x_1^2 + 2x_1x_2} P - \sigma \leq 0, g_3 = \frac{1}{x_1 + \sqrt{2}x_2} P - \sigma \leq 0,
\]

where \(0 \leq x_1 \leq 1\), \(0 \leq x_2 \leq 1\), \(H = 100\,\text{cm}\), \(P = 2\,\text{KN/cm}^2\), and \(\sigma = 2\,\text{KN/cm}^2\).

\begin{figure}[H] 
	\centering 
	\includegraphics[width=0.4\textwidth]{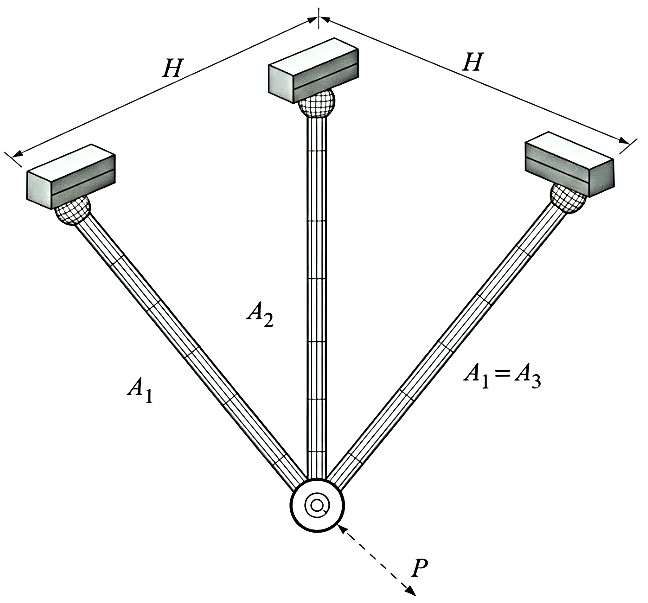} 
	\renewcommand{\figurename}{Figure}
	\caption{Schematic diagram of the three-bar truss design (TTD) problem.} 
	\label{Fig.TTD} 
\end{figure}

\begin{table}[htbp]
	\centering
	\renewcommand{\tablename}{Table}
	\caption{Results analysis of VDO on TTD.}
	\label{TTD}
	{\scriptsize
		\begin{tabular}{llllll}
			\toprule
			\textbf{Algorithm} & \textbf{x1} & \textbf{x2} & \textbf{Best} & \textbf{Mean} & \textbf{Variance} \\
			\midrule
			VDO   & 7.887E-01 & 4.082E-01 & 2.638958E+02 & 2.638958E+02 & \textbf{1.3699E-20} \\
			GA    & 7.891E-01 & 4.071E-01 & 2.638961E+02 & 2.639055E+02 & 5.7860E-05 \\
			PSO   & 7.887E-01 & 4.082E-01 & 2.638958E+02 & 2.638959E+02 & 1.2459E-09 \\
			GWO   & 7.883E-01 & 4.092E-01 & 2.638959E+02 & 2.638970E+02 & 6.4164E-07 \\
			WOA   & 7.878E-01 & 4.108E-01 & 2.638965E+02 & 2.641168E+02 & 7.9991E-02 \\
			HHO   & 7.884E-01 & 4.089E-01 & 2.638959E+02 & 2.639401E+02 & 3.9922E-03 \\
			HO    & -2.9540E+00 & -2.0546E+00 & -1.04099E+03 & -4.785923E+02 & 1.6850E+05 \\
			NRBO  & 7.887E-01 & 4.082E-01 & 2.638958E+02 & 2.638958E+02 & 3.2679E-25 \\
			FATA  & 7.874E-01 & 4.120E-01 & 2.638974E+02 & 2.639067E+02 & 6.3384E-05 \\
			MGO   & 7.887E-01 & 4.083E-01 & 2.638958E+02 & 2.638976E+02 & 1.7484E-05 \\
			LEA   & 7.855E-01 & 4.175E-01 & 2.639077E+02 & 2.641587E+02 & 6.5939E-02 \\
			ALA   & 7.887E-01 & 4.082E-01 & 2.638958E+02 & 2.638958E+02 & \textbf{7.9108E-27} \\
			\bottomrule
	\end{tabular}}
\end{table}

Table~\ref{TTD} reportes the optimization performance of VDO and eleven comparative algorithms on the three-bar truss design (TTD) problem. The proposed Virus Diffusion Optimizer (VDO) achieves the best overall performance, obtaining the lowest mean objective value of $2.638958 \times 10^{2}$ and the smaller variance of $1.37 \times 10^{-20}$, demonstrating exceptional convergence precision and remarkable numerical stability. Classical methods such as GA and PSO yield slightly higher mean values and noticeably higher variances ($5.79 \times 10^{-5}$ and $1.25 \times 10^{-9}$, respectively), suggesting that they are more prone to small-scale oscillations during convergence. Advanced metaheuristics like GWO, HHO, and MGO also perform well in terms of accuracy but still fail to match VDO’s consistency. The HO algorithm notably diverge, producing negative infeasible results due to instability in constraint handling.

Overall, the results confirm that VDO exhibits outstanding reliability and precision in this problem. Its virus-inspired dual-diffusion mechanism effectively enhance both exploitation accuracy and search stability, enabling the algorithm to consistently achieve optimal feasible designs with minimal variation across independent runs.

\subsection{Welded-beam design}
The welded beam is a common engineering optimisation problem with an objective to find an optimal set of the dimensions 
Decision variables:

(i) Shell thickness ($h = x_1$)

(ii) Head thickness ($l = x_2$)

(iii) Inner radius ($t = x_3$)

(iv) Length of the vessel without including the head ($b = x_4$)

such that the fabrication cost of the beam is minimized. This problem is a continuous optimisation problem. See the Figure below for graphical details of the beam dimensions (h,l,t,b) to be optimised.

Minimize:
\[f(x) = 1.10471{x_1^2}{x_2} + 0.04811{x_3}{x_4}(14 + x_2) \]

Subject to:
\[ g_1(\vec{x}) = \tau(\vec{x}) - \tau_{\max} \leq 0, \  g_2(\vec{x}) = \sigma(\vec{x}) - \sigma_{\max} \leq 0, \ 
g_3(\vec{x}) = x_1 - x_4 \leq 0, \]
\[ g_4(\vec{x}) = 0.10471x_1^2 + 0.04811x_3x_4(14 + x_2) - 5 \leq 0, \]
\[ g_5(\vec{x}) = 0.125 - x_1 \leq 0, \ 
g_6(\vec{x}) = \delta(\vec{x}) - \delta_{\max} \leq 0, \  g_7(\vec{x}) = P - P_c(\vec{x}) \leq 0, \]
With bounds:
\[ 0.1 \leq x_1 \leq 2, \  0.1 \leq x_2 \leq 10, 0.1 \leq x_3 \leq 10, \ 0.1 \leq x_4 \leq 2. \]

where:

\[ \tau(\vec{x}) = \sqrt{(\tau')^2 + 2\tau' \tau'' \frac{x_2}{2R} + (\tau'')^2}, \tau' = \frac{P}{\sqrt{2x_1 x_2}}, \quad \tau'' = \frac{MR}{J}, \quad M = P\left(L + \frac{x_2}{2}\right), \]
\[ R = \sqrt{\frac{x_2^2}{4} + \frac{(x_1 + x_3)^2}{4}}, J = 2 \left[ \sqrt{2x_1 x_2} \left( \frac{x_2^2}{12} + \frac{(x_1 + x_3)^2}{4} \right) \right], \]
\[ \sigma(\vec{x}) = \frac{6PL}{x_4 x_3^2}, \ \delta(\vec{x}) = \frac{4PL^3}{E x_3^3 x_4},  P_c(\vec{x}) = \frac{4.013 E \sqrt{\dfrac{x_3^2 x_4^6}{36}}}{L^2} \left( 1 - \frac{x_3}{2L} \sqrt{\frac{E}{4G}} \right),\]
\[ P = 6000 \, \text{lb}, \quad L = 14 \, \text{in}, \quad E = 30 \times 10^6 \, \text{psi}, G = 12 \times 10^6 \, \text{psi}, \]
\[ \tau_{\max} = 13{,}600 \, \text{psi}, \quad \sigma_{\max} = 30{,}000 \, \text{psi}, \quad \delta_{\max} = 0.25 \, \text{in}. \]

\begin{figure}[H] 
	\centering 
	\includegraphics[width=0.4\textwidth]{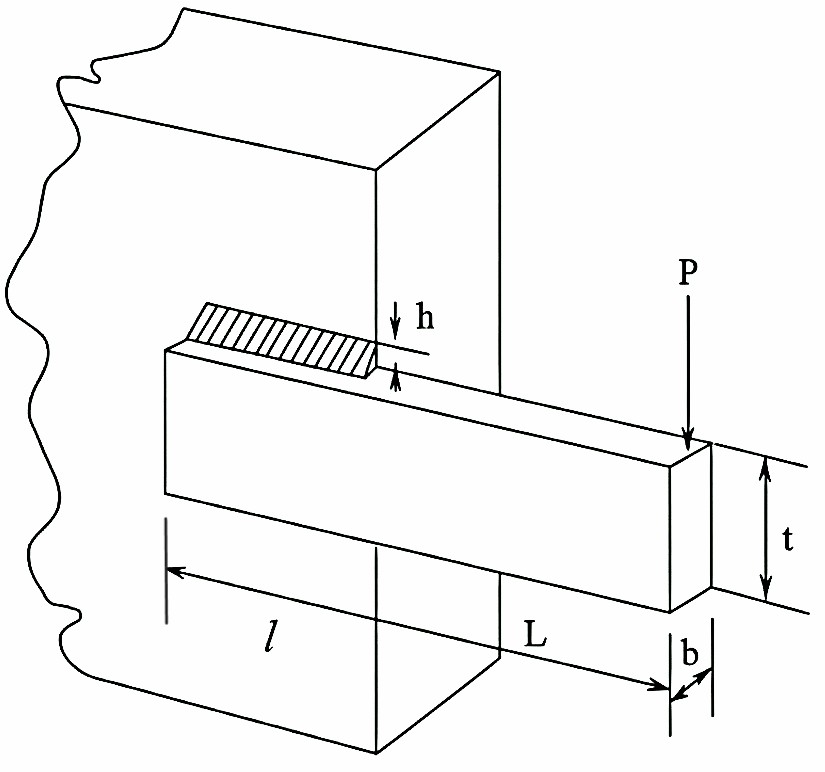} 
	\renewcommand{\figurename}{Figure}
	\caption{Schematic diagram of the welded beam design (WBD) problem.} 
	\label{Fig.WBD} 
\end{figure}

\begin{table}[htbp]
	\centering
	\renewcommand{\tablename}{Table}
	\caption{Results analysis of VDO on WBD.}
	\label{WBD}
	{\scriptsize
		\begin{tabular}{llllllll}
			\toprule
			\textbf{Algorithm} & \textbf{x1} & \textbf{x2} & \textbf{x3} & \textbf{x4} & \textbf{Best} & \textbf{Mean} & \textbf{Std} \\
			\midrule
			VDO   & 2.0570E-01 & 3.2349E+00 & 9.0366E+00 & 2.0570E-01 & \textbf{1.6928E+00} & \textbf{1.6928E+00} & \textbf{1.2082E-27} \\
			GA    & 2.0680E-01 & 3.2243E+00 & 9.0128E+00 & 2.0690E-01 & 1.6972E+00 & 1.9583E+00 & 8.1948E-02 \\
			PSO   & 2.0570E-01 & 3.2349E+00 & 9.0366E+00 & 2.0570E-01 & \textbf{1.6928E+00} & \textbf{1.6928E+00} & 3.4303E-14 \\
			GWO   & 2.0560E-01 & 3.2372E+00 & 9.0372E+00 & 2.0570E-01 & 1.6930E+00 & 1.6937E+00 & 2.9533E-07 \\
			WOA   & 1.9290E-01 & 3.5382E+00 & 9.0027E+00 & 2.0990E-01 & 1.7401E+00 & 2.0443E+00 & 2.5620E-01 \\
			HHO   & 2.0760E-01 & 3.2439E+00 & 8.9955E+00 & 2.0760E-01 & 1.7038E+00 & 1.7962E+00 & 4.6702E-03 \\
			HO    & 2.0370E-01 & 3.2780E+00 & 9.0366E+00 & 2.0570E-01 & 1.6956E+00 & 1.8188E+00 & 3.0930E-02 \\
			NRBO  & 2.0580E-01 & 3.2337E+00 & 9.0342E+00 & 2.0580E-01 & 1.6932E+00 & 1.7325E+00 & 8.4343E-04 \\
			FATA  & 2.0570E-01 & 3.2357E+00 & 9.0366E+00 & 2.0570E-01 & 1.6929E+00 & 1.6939E+00 & 1.5248E-06 \\
			MGO   & 2.0570E-01 & 3.2352E+00 & 9.0366E+00 & 2.0570E-01 & \textbf{1.6928E+00} & 1.6931E+00 & 1.8126E-07 \\
			LEA   & 2.0680E-01 & 3.3470E+00 & 8.7542E+00 & 2.1970E-01 & 1.7634E+00 & 1.8991E+00 & 1.7145E-02 \\
			ALA   & 2.0570E-01 & 3.2349E+00 & 9.0366E+00 & 2.0570E-01 & \textbf{1.6928E+00} & \textbf{1.6928E+00} & \textbf{9.5377E-31} \\
			\bottomrule
	\end{tabular}}
\end{table}

The welded beam design (WBD) problem is a standard constrained continuous optimisation benchmark that aims to minimise the fabrication cost of a welded beam subject to strength, deflection, and buckling constraints. Table~\ref{WBD} present the optimisation results obtained by the proposed Virus Diffusion Optimizer (VDO) and eleven comparative algorithms. The results show that VDO achieves the global optimum design parameters $x_1 = 0.2057$, $x_2 = 3.2349$, $x_3 = 9.0366$, and $x_4 = 0.2057$, leading to the minimum fabrication cost of ${1.6928}$ with an extremely small variance of ${1.21\times10^{-27}}$. These results precisely satisfy the design constraints and represent a physically feasible and cost-optimal configuration.

Compared with the benchmark algorithms, VDO demonstrate clear superiority in both accuracy and consistency. While several algorithms such as PSO and GWO also approach the same best value ($1.6928$), their reported variances ($3.43\times10^{-14}$ and $2.95\times10^{-7}$, respectively) are either higher or less stable across runs. In contrast, VDO consistently converge to the optimum solution in all independent runs, achieving near-zero standard deviation and thereby confirming its exceptional robustness and reproducibility. Classical algorithms such as GA and WOA produce noticeably worse results, with mean costs of $1.9583$ and $2.0443$, respectively, and standard deviations larger by several orders of magnitude, highlighting their sensitivity to local minima and limited convergence precision.

Overall, the WBD results substantiate that VDO provides highly accurate, stable, and reproducible optimisation performance. Its virus-inspired search dynamics allow it to maintain a delicate balance between exploration and exploitation, efficiently navigating the nonlinear constrained landscape of the welded beam design problem. Consequently, VDO not only achieve the global optimum but also deliver better numerical precision, demonstrating its strong potential for real-world engineering design optimisation tasks where solution stability and constraint satisfaction are critical.

\section{Conclusion}

This paper proposed the Virus Diffusion Optimizer (VDO), a novel metaheuristic algorithm inspired by the lifecycle of the Herpes Simplex Virus (HSV). By emulating viral tropism, replication regulation, virion diffusion, and latency reactivation, VDO achieved an effective balance between exploration and exploitation in high-dimensional search spaces. 
Extensive evaluations on CEC2017, CEC2022 benchmark functions and engineering design problems—including pressure vessel design (PVD), three-bar truss design (TTD), and welded beam design (WBD)—demonstrated VDO's superior performance over 11 algorithms. The key strengths of VDO could be summarized as follows:

\begin{itemize}
	\item \textbf{Strong Robustness:} The algorithm exhibited low standard deviation across runs, highlighting its reliability.
	\item \textbf{Effective Search Balance:} Adaptive mechanisms dynamically maintained a harmony between global exploration and local exploitation.
	\item \textbf{Sustained Population Diversity:} The latency reactivation strategy mitigated premature convergence and preserved search potential.
\end{itemize}

However, VDO also presented certain limitations. The algorithm's implementation was relatively complex, resulting in higher computational costs, particularly as problem scale increases. Future work focused on simplifying hyperparameters without compromising performance, developing hybrid strategies with local search techniques, and extending VDO to multi-objective optimization domains. In conclusion, VDO offered a robust and efficient bio-inspired framework for complex optimization tasks, underscoring the continued relevance of biological metaphors in the design of metaheuristic algorithms.

\section*{Declarations}
{\bf Funding:} This research is supported by the National Natural Science Foundation of China (No. 12401415), the 111 Project (No. D23017), the Natural Science Foundation of Hunan Province (No. 2025JJ60009). 

\noindent{\bf Data Availability:} Enquiries about code availability should be directed to the authors.

\noindent{\bf Competing interests:} The authors have no competing interests to declare that are relevant to the content of this article.

\bibliographystyle{abbrv}
\bibliography{refs}

\end{document}